\documentclass[a4paper]{article}

\usepackage[utf8]{inputenc}
\usepackage[english]{babel}
\usepackage[T1]{fontenc}
\usepackage{amsmath}
\usepackage{mathrsfs}
\usepackage{amsfonts}
\usepackage{amssymb}
\usepackage{graphicx}
\usepackage{enumitem}
\usepackage{fancyhdr}
\usepackage{listings}
\usepackage[symbol]{footmisc}
\usepackage{hyperref}
\hypersetup{
    colorlinks=true,
    linkcolor=red,
    filecolor=magenta,      
    urlcolor=cyan,
    linktoc=page
    }
\usepackage{url}
\usepackage{mathabx}
\usepackage{tikz-cd}
\usepackage{tkz-euclide}
\usetikzlibrary{math}
\usetikzlibrary{decorations.pathreplacing,angles,quotes}
\usetikzlibrary{arrows}
\usetikzlibrary{braids}
\usepackage{pgfplots}
\pgfplotsset{compat=1.15}
\usepackage{mathrsfs}

\usepackage{footnote}       %to have the footnote in the tabular
\makesavenoteenv{tabular}   %to have the footnote in the tabular
    
\usepackage[backend=biber, maxnames=1000, minnames=100]{biblatex}
\usepackage{amsthm}
\usepackage{perpage}
\usepackage{multicol}
\usepackage{multirow}
\usepackage[justification=centering]{caption}
\usepackage{wrapfig}
\usepackage[normalem]{ulem}
\usepackage{listings}
\usepackage{comment}

\usepackage{pdflscape}
\usepackage{afterpage}
\usepackage{capt-of}
\usepackage{bm}
\usepackage{csquotes}

\newtheoremstyle{bfnote}%
{}{}%
{\itshape}{}%
{\bfseries}{.}%
{ }%
{\thmname{#1}\thmnumber{ #2}\thmnote{ -- #3}}

\theoremstyle{bfnote}

\newtheorem{theorem}{Theorem}[section]

\newtheorem{proposition}[theorem]{Proposition}
\newtheorem{lemma}[theorem]{Lemma}
\newtheorem{coro}[theorem]{Corollary}
\newtheorem{conj}[theorem]{Conjecture}

\newtheorem{definition}[theorem]{Definition}
\newtheorem{def/prop}[theorem]{Definition/Proposition}

\newtheoremstyle{bfnote_def}%
{}{}%
{}{}%
{\bfseries}{.}%
{ }%
{\thmname{#1}\thmnumber{ #2}\thmnote{ -- #3}}

\theoremstyle{bfnote_def}

\newtheorem{remark}{Remark}
\newtheorem{openquestion}{Open question}

\newcommand{\Zn}{\mathbb{Z}/n\mathbb{Z}}

\newcommand{\PGL}[1]{\mathbb{P}\mathrm{GL}_{#1}}
\newcommand{\Tn}{\mathcal{T}_{n}}
\newcommand{\tildeTn}{\tilde{\mathcal{T}}_{n}}

\newcommand{\tildePn}{\tilde{\mathcal{P}}_{n}}
\newcommand{\floor}[1]{\left\lfloor #1 \right\rfloor}

\title{Collapsing in polygonal dynamics}

\author{Jean-Baptiste Stiegler\footnote{Laboratoire de mathématiques d'Orsay, Équipe Topologie et Dynamique, 91405, Orsay, France. \textit{E-mail address}: \texttt{jean-baptiste.stiegler at universite-paris-saclay.fr}}}
\date{\today}

\begin{document}

\maketitle

\abstract{We define polygonal dynamics as a family of dynamical systems acting on points in projective spaces. The most famous example is the pentagram map. Similar collapsing phenomena seem to occur in most of these systems. We prove it in some case, and conjecture that it almost always happens. Moreover, we give a formula for the limit point in term of roots of $d+1$ degree polynomials (where $d$ is the dimension of the projective space). We do so by generalizing Glick's operator, interpreted as an infinitesimal monodromy. This answers questions about its reappearance in many systems, together with preserved quantities. We apply these results to several polygonal dynamics in $\mathbb{P}^1$ and introduce a new one called ``staircase'' cross-ratio dynamics, for which we study particular configurations.}

\tableofcontents

\section{Introduction}

Before starting, a word on vocabulary. We denote the projective spaces by $\mathbb{P}^d$ and their automorphism groups by $\PGL{d+1}$ without specifying over which field $F$ they are defined. This is because we see them from the algebraic geometry standpoint, as schemes (like it's done in \cite{weinreichAlgebraicDynamicsPentagram2023}). Similarly, the group scheme $\mathbb{G}_m$ refers to the group of invertible elements $F^*$. However, no algebraic geometry background will be needed here, and the unfamiliar reader may consider the statements over a fixed field, like $\mathbb{C}$.

\subsection{Definitions}
The point of this paper is to study a family of  dynamical systems, called \textit{polygonal dynamics}. It stems from observing similar behaviours when numerically simulating them. Let's start with a definition.

\begin{definition}\label{def:poly_dyn}
    Let $n,d$ be positive integers. A \textbf{twisted $\bm{n}$-gon} $P$ is the data of \textbf{vertices} $p_1,\dots,p_n \in \mathbb{P}^d$ and a \textbf{monodromy} $M\in\PGL{d+1}$. The vertices are to be thought as an infinite sequence satisfying:
    \[p_{i+n}=Mp_i\quad \forall i\in\mathbb{Z}.\]
    Moreover, we impose a \textbf{non-degeneracy condition}, meaning that for any $i\in\mathbb{Z}$, the vertices $(p_i,\dots,p_{i+d+1})$ form a projective frame.
    
    $P$ is said be \textbf{closed} if $M=\mathrm{Id}$.  The set of twisted and closed $n$-gons are respectively denoted by $\tildeTn$ and $\tildePn$.

    A \textbf{polygonal transformation} is a birational map \[T:\tildeTn\dashrightarrow \tildeTn\]
    which commutes with the following action of $\PGL{d+1}$:
    \[A\cdot P=((A p_1,\dots,Ap_n), AMA^{-1}),\quad A\in \PGL{d+1}.\]
    The class of $P$ is denoted by $[P]$. The action of $\PGL{d+1}$ allows to define the quotient spaces $\mathcal{T}_n$ and $\mathcal{P}_n$, also called moduli spaces.
\end{definition}

Since $T$ commutes with the action of $\PGL{d+1}$, it also yields a birational map (again denoted by $T$) over the moduli spaces. In many cases, the dynamic of $T$ on the moduli spaces have a very strong property called discrete integrability.\footnote{The following definition works for the algebraic case, and is inspired from \cite{weinreichAlgebraicDynamicsPentagram2023}. For general review of continuous integrable systems, and the link between the Hamiltonian and algebraic flavors, see \cite{audin:hal-00135277}.}

\begin{definition}\label{def:integr}
    A dynamical system $(X,T)$ is said to be \textbf{discretely integrable} if there is a family $A$ of abelian varieties\footnote{Here, this is the algebraic counterpart of lagrangian tori.} which foliate $X$ such that, on a dense subset $\Sigma \subset X$, there is a birational map $\delta:\Sigma \to A$ which identifies some iterate of $T$ with a relative translation $\tau: A\to A$.
\end{definition}

\begin{remark}
    Definition~\ref{def:poly_dyn} can be generalized by replacing $\mathbb{P}^d$ with a topological space $X$ and $\PGL{d+1}$ by a group $G$ acting by homeomorphism (see \cite[§2]{beffaIntegrableGeneralizationsPentagram2015}). For instance, the polygon recutting defined by Adler \cite{adlerRecuttingsPolygons1993} in 1993 takes place in $X=\mathbb{R}^2$ with $G=\mathrm{Isom}^+(\mathbb{R}^2)$. See table~\ref{tab:examples_polygonal_dynamics} for a review of existing polygonal dynamics. Since our results only concerns polygonal dynamics in $\mathbb{P}^d$, we will set the rest aside.
\end{remark}

\subsection{Previous works}

\begin{figure}
    \centering
    \includegraphics[width=0.5\linewidth]{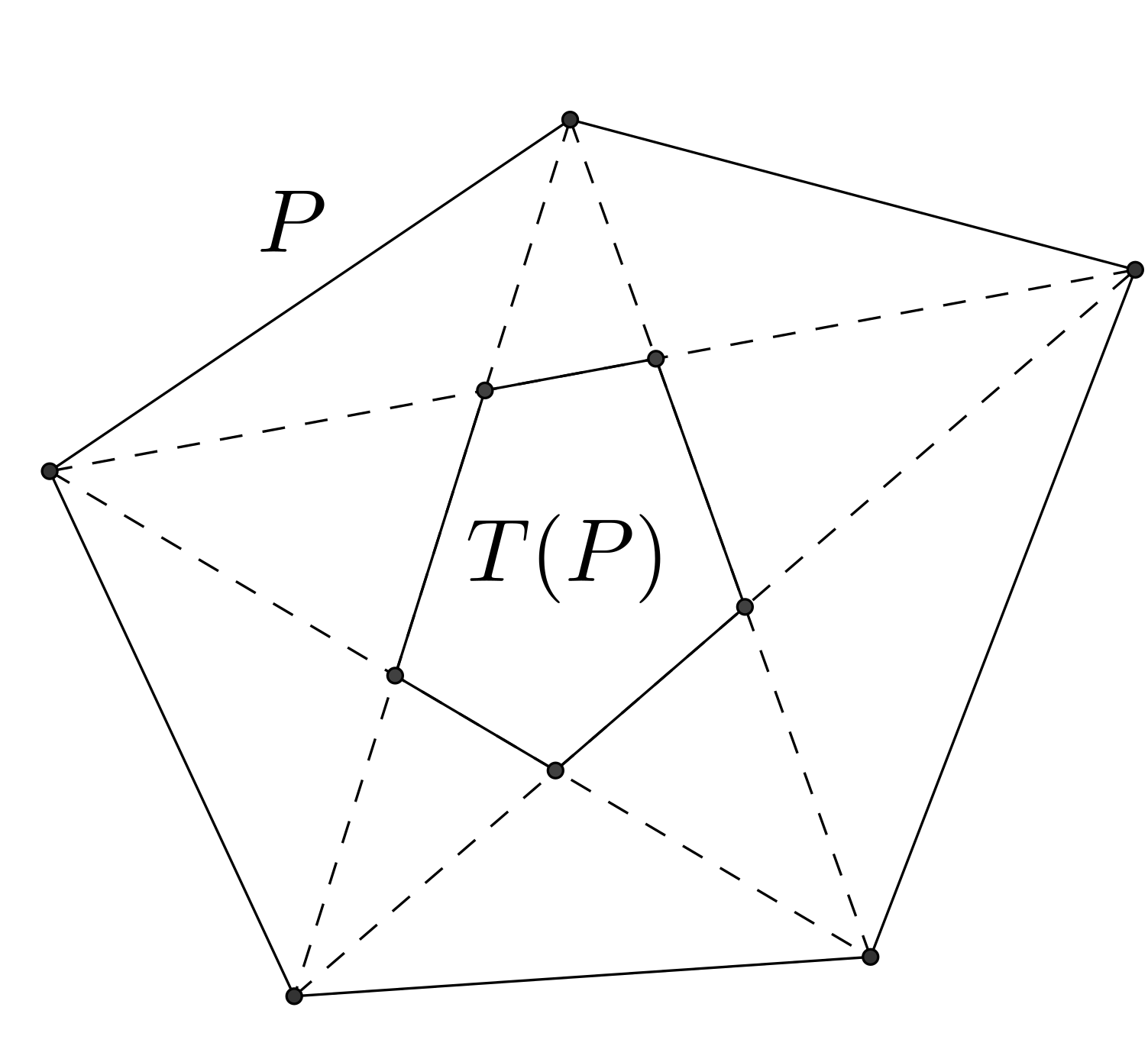}
    \caption{The pentagram map applied on a closed pentagon in $\mathbb{P}^2(\mathbb{R})$.}
    \label{fig:pentagram}
\end{figure}

Let's start with the most famous polygonal dynamic: the \textit{pentagram map}, defined by Schwartz \cite{schwartzPentagramMap1992} in 1992. In his original paper, he defines the map $T$ acting on the set of strictly convex closed $n$-gons (with $n\geq5$, hence the name) in $\mathbb{P}^2(\mathbb{R})$.\footnote{He generalised it to twisted polygons \cite{schwartzDiscreteMonodromyPentagrams2007} in 2007.} A polygon $P=(p_1,\dots,p_n)$ is sent to another one by taking the intersections of the small diagonals $\overline{p_{i-1} p_{i+1}}$ and $\overline{p_{i} p_{i+2}}$, where the indices are understood modulo $n$ (see figure~\ref{fig:pentagram} for an example). Schwartz proved that for strictly convex polygons, the sequence of iterates collapses exponentially fast towards a point. However, the meaning of this point, or even the existence of a formula to express it from the original vertices of $P$, remained unknown for 26 years. This collapse point, which can be thought as a kind of ``projective center of mass'' of the polygon $P$, is quite intriguing.

In 2018, Glick \cite{glickLimitPointPentagram2020} gives a formula for the collapse point. To do this, he defines a linear operator of $\mathbb{R}^3$ associated to a convex polygon $P$, now called \textit{Glick's operator}, which is invariant by the dynamic and sends $P$ to its interior. He concludes that the collapse point must be a fixed point of the operator. The formula is quite simple, involving polynomials of degree~$3$ whose coefficients depend on the vertices of $P$. However, the proof is quite computational and doesn't shed much light on the meaning of the point. 

In 2020, Aboud and Izosimov \cite{aboudLimitPointPentagram2022} reinterpret Glick's operator in term of a matrix $M'$, called an \textit{infinitesimal monodromy}. To do this, they ``open'' the closed polygon into a family of twisted ones through a tool called a \textit{scaling symmetry}.

\begin{definition}\label{def:scaling_sym}
A \textbf{scaling symmetry} is a non-trivial group action
\[\rho:\mathbb{G}_m\times\Tn\to\Tn\]
which commutes with the dynamic of $T$.
\end{definition}

In 2023, Arnold and Arreche \cite[§7.2]{arnoldSymmediansHyperbolicBarycenters2024} notice a similarity between the infinitesimal monodromy of the pentagram map and a so-called ``hyperbolic barycenter'' defined for the flat cross-ratio dynamics (see \cite[§6.2.1]{arnoldCrossratioDynamicsIdeal2022}). 

In 2024, Schwartz \cite[§9.2]{schwartzFlappingBirdsPentagram2026} investigates a family of transformations $\Delta_k$ (the pentagram map being $\Delta_1$) for which certain polygons collapse. Moreover, he conjectures that the collapsing point is given by a cousin of Glick's operator. In the first version of his preprint\footnote{\href{2403.05735v1}{arXiv:2403.05735v1} page 57.}, he makes the following observation: 
\begin{displayquote}
    ``\textit{I wonder if this means that the collapse point exists for all starting points of the pentagram map. Even if the iterations go completely crazy under the map, perhaps they still collapse to the point predicted by Glick’s operator. The idea of a completely general collapse point has always seemed absurd to me, but maybe it is not. Nobody knows.}''
\end{displayquote}

\subsection{Plan of the paper and main results}

In this paper, we aim to answer these questions by considering the general setting of polygonal dynamics in $\mathbb{P}^d$, for any $d\in\mathbb{N}^*$. Examples of such systems are presented in Table \ref{tab:examples_polygonal_dynamics}, located in the appendix.

In Section \ref{section:def}, we construct a coordinate system on the moduli space using a generalization of the usual cross-ratio. This allows us to determine the monodromy of any polygon. Then, if the system admits a scaling symmetry, we use it to deform polygons. In particular, we generalize the computation of \cite{aboudLimitPointPentagram2022} and give an explicit formula of the infinitesimal monodromy in term of the cross-ratio coordinates (Theorem \ref{thm:infinitesimal_monod}). If the scaling symmetry acts linearly,\footnote{Meaning that the cross-ratio coordinates vary linearly, which is the case in every examples.} we compute higher derivatives of the monodromy.

In Section \ref{section:in P^n}, we investigate the collapsing of polygons, meaning that all the vertices converge towards the same point. Our computer simulations seem to indicate that it almost always happens when the dynamic takes place over $\mathbb{C}$. We manage to prove it for polygons whose orbit is periodic on the moduli space (Theorem \ref{thm:asympt_periodic}), which are known to exist. Moreover, we state in Conjecture \ref{conj:asympt} that it is the case for almost every polygons. In Theorem \ref{thm:collapse_closed}, we prove that the collapse point of a closed polygon has to be a fixed point of the infinitesimal monodromy, which answers the question of Schwartz in \cite[§9.2]{schwartzFlappingBirdsPentagram2026}. We obtain in Corollary \ref{cor:field_ext} that the coordinates of the limit point are roots of degree $d+1$ polynomials whose coefficients are determined by the polygon $P$.

In Section \ref{section:in P^1}, we focus on dynamics in $\mathbb{P}^1$ and give effective versions of the previous results. In Proposition \ref{prop:act_trivial}, we prove that any dynamic is trivial on the moduli spaces of closed triangles and quadrilaterals. In Theorem \ref{thm:monod_closed}, we give a more explicit formula of the infinitesimal monodromy, which is exactly the matrix which appeared in \cite[§6.2.1]{arnoldCrossratioDynamicsIdeal2022}). This answers to the question of Arnold and Arreche \cite[§7.2]{arnoldSymmediansHyperbolicBarycenters2024}. We also write explicitly the quadratic polynomial for which the limit point vanishes. We conclude the section by deriving the formulas for two sets of $\floor{n/2}$ preserved quantities (Propositions \ref{prop:Hk} and \ref{prop:Gk}), which are good candidates for Hamiltonians.

In Section \ref{section:examples}, we apply these results to two polygonal dynamics in $\mathbb{P}^1$, namely the leapfrog map (defined in \cite[§5.2]{gekhtmanIntegrableClusterDynamics2016}) and the ``flat'' cross-ratio dynamic (studied in \cite{hetrich-jerominPeriodicDiscreteConformal2001},  \cite{arnoldCrossratioDynamicsIdeal2022} and \cite{affolterIntegrableDynamicsProjective2023}). We find their scaling symmetries (Propositions \ref{prop:leapfrog_scaling} and \ref{prop:flat_scaling}) and derive their associated infinitesimal monodromy and limit point.

In Section \ref{section:staircase}, we introduce a new polygonal dynamic on $\mathbb{P}^1$, called the ``staircase'' cross-ratio dynamics. As the name shows, it is a cousin of the ``flat'' cross-ratio dynamics. They both stem from the theory of discrete holomorphic functions and their continuation (see \cite[§6]{affolterSchwarzianOctahedronRecurrence2024} and references therein). They differ by a choice of geometry (flat or staircase like) for the initial conditions. The ``staircase'' cross-ratio dynamics consist of a group of polygonal transformations called ``flips'', which happens to be the affine symmetric group (Theorem \ref{thm:affine_sym_grp}). As for the other dynamics, we find a scaling symmetry and derive its infinitesimal monodromy and limit point. We conclude with a family of $n$-gons with periodic orbits on their moduli space, whose vertices always collapse toward the fixed points of the infinitesimal monodromy, as predicted by Theorem \ref{thm:collapse_closed}.

\section{Basic properties}\label{section:def}

We'll make an extensive use of the first fundamental theorem of projective geometry (see \cite[Prop. 4.5.10]{BergerMarcel1927-20162016GT1}), which states that given two projective frames, there exists a unique projective transformation mapping the first to the second.

\begin{lemma}
    Polygonal transformations preserve monodromies.
\end{lemma}

\begin{proof}
    Let $P=((p_i)_{i\in\mathbb{Z}},M)$ be an $n$-gon in $\mathbb{P}^d$, $T$ be a polygonal transformation, and $T(P)=:P'=((p'_i)_{i\in\mathbb{Z}},M')$. We have 
    \[M\cdot P'=((Mp'_i)_{i\in\mathbb{Z}},MM'M^{-1}),\]
    and by the monodromy condition
    \[M\cdot P = ((p_{i+n})_{i\in\mathbb{Z}},M).\]
    Since $T$ commutes with the action of $\PGL{d+1}$
    \[M\cdot P'=T(M\cdot P)=((p'_{i+n})_{i\in\mathbb{Z}},M'),\]
    so we get that for all $i\in\mathbb{Z}$
    \[M'p'_i=p'_{i+n}=Mp'_i.\]
    Since $M$ and $M'$ map the projective frame $(p'_0,\dots,p'_{d+1})$ to $(p'_n,\dots,p'_{d+1+n})$, we conclude by the first fundamental theorem of projective geometry that they are equal.
\end{proof}

\subsection{Coordinate system in the moduli space}
First of all, we need a way to parameterize the moduli space. We do so by generalizing the usual cross-ratio in $\mathbb{P}^1$ to $\mathbb{P}^d$.

\begin{definition}\label{def:cr}
    Let $p_0,\dots,p_{d+1} \in \mathbb{P}^d$ be points forming a projective frame. There exists a unique projective transformation $C_{p_0,\dots,p_{d+1}}$ which maps $(p_0,p_1,\dots,p_{d+1})$ to the standard frame
    \[([1:\dots:1],[1:0:\dots:0],\dots,[0:\dots:0:1]).\]
    If $p$ is another point in $\mathbb{P}^d$, then its \textbf{cross-ratio with respect to the frame} $\bm{(p_0,...,p_{d+1})}$ is \[c=[x_1:\dots:x_{d+1}]:=C_{p_0,\dots,p_{d+1}}(p).\]
\end{definition}

The following proposition derives from the definition, and is proved exactly in the same way as for the $\mathbb{P}^1$ case (see \cite[Prop. 6.2]{AudinMichèle2006G} for details).

\begin{proposition}\label{prop:cr_inv}
    This cross-ratio is projectively invariant. Moreover, there exist a projective transformation mapping a frame $(p_0,\dots,p_{d+1},p)$ to another one $(p_0',\dots,p_{d+1}',p')$ if and only if their respective cross-ratios coincide.
\end{proposition}

There are explicit formulas involving determinants, using the so called volume cross-ratio defined in \cite{Olver_2001}. For a geometric understanding involving only the usual cross-ratio on a projective line, see \cite[Prop. 6.5.11]{BergerMarcel1927-20162016GT1}.  

\begin{proposition}\label{prop:volume_cr}
    Let $p_0,\dots,p_{d+1}$ be a projective frame in $\mathbb{P}^d$, $p\in \mathbb{P}^d$ be another point, and  $V_0,\dots,V_{d+1},V$ be a lift in $\mathbb{A}^{d+1}$. Let $[x_1:\dots:x_{d+1}]$ be the cross ratio defined previously.
    For any $1\leq i<j\leq d+1$, call $V^{i,j}$ the collection of all vectors from $V_1$ to $V_{d+1}$, excluding $V_i$ and $V_j$. Then we have: 
    \[\frac{x_i}{x_j}=\frac{\mathrm{det}(V^{i,j},V_i,V_0)}{\mathrm{det}(V^{i,j},V_i,V)}\frac{\mathrm{det}(V^{i,j},V_j,V)}{\mathrm{det}(V^{i,j},V_j,V_0)}.\]
\end{proposition}
\begin{remark}
    This retrieves the usual definition of the cross-ratio of four points in $\mathbb{P}^1$:
    \[[p_1,p_2,p_0,p]=\frac{p_1-p_0}{p_1-p}\frac{p_2-p}{p_2-p_0}.\]
\end{remark}

\begin{proof}
    First of all, we see that any choice of lift gives the same quantity. Indeed, if any vector is multiplied by a non-zero scalar, it cancels out in the ratio of determinants. Then, we see it is projective invariant, as if we apply any linear automorphism $A$ to the vectors, all the determinants get multiplied by $\mathrm{det}(A)$, which in turn cancels out in the ratio. Hence, we can work with
    \[V_0=\begin{pmatrix}
       1\\
       \vdots\\
       1
    \end{pmatrix},\quad V_1=\begin{pmatrix}
       1\\
       0\\
       \vdots\\
       0
    \end{pmatrix},\quad\dots,\quad V_{d+1}=\begin{pmatrix}
       0\\
       \vdots\\
       0\\
       1
    \end{pmatrix},\quad V=\begin{pmatrix}
       x_1\\
       \vdots\\
       x_{d+1}
    \end{pmatrix}.\]
    We conclude by the straight out computation:
    \[\frac{\mathrm{det}(V^{i,j},V_i,V_0)}{\mathrm{det}(V^{i,j},V_i,V)}\frac{\mathrm{det}(V^{i,j},V_j,V)}{\mathrm{det}(V^{i,j},V_j,V_0)}=\frac{1\times x_i}{x_j\times 1}.\qedhere\]
\end{proof}

\begin{definition}\label{def:cr_coord}
Let $P$ be an $n$-gon in $\mathbb{P}^d$. Its \textbf{cross-ratio coordinates} are defined to be $C=(c_i)_{i\in\mathbb{Z}}$, where
\[c_i=[c_{i,1}:\dots:c_{i,d+1}]:=C_{p_{i-1},p_{i},\dots,p_{i+d}}(p_{i+d+1}).\]
They are $n$-periodic, and because of the non-degeneracy each $c_{i,j}$ is non-zero. We can normalize each $c_{i,d+1}$ to $1$, and the formula from proposition~\ref{prop:volume_cr} gives us the values of the others $c_{i,j}\in \mathbb{P}^1\setminus\{0,\infty\}$. 
\end{definition}

\begin{remark}
    Other coordinate systems can be used to study polygonal dynamics, for instance the corner coordinates in the case of the pentagram map (see \cite[p.415]{ovsienkoPentagramMapDiscrete2010}). However, \cite[Thm. 3.10]{Olver_2001} tells us that any coordinate system is equivalent to the one given by cross-ratios. We will make concrete translation between different coordinates systems in section \ref{section:examples}.
\end{remark}

The cross-ratio coordinates $C$ parameterise the projective class $[P]$ of a non-degenerate polygon. With the following lemma, which is a generalization of \cite[Lemma 3.2]{arnoldCrossratioDynamicsIdeal2022}, we can in fact reconstruct a polygon from the cross-ratio coordinates. Hence $\Tn$ is isomorphic to $(\mathbb{P}^1\setminus\{0,\infty\})^{dn}$. For the quotient space including degenerate polygons, see \cite{weinreichGITStabilityLinear2024a}.

\begin{lemma}\label{lemma:monodromy}
    Let $P$ be a $n$-gon in $\mathbb{P}^d$ and
    \[c_{i}=[c_{i,1}:\dots:c_{i,d+1}]\]
    be its cross-ratio coordinates. For any $i\in\{1,\dots,n\}$, define the matrix
    \[L_i=\begin{pmatrix}
       0 & \dots & 0 & c_{i,1}\\
       -c_{i,2} & \dots & 0 & c_{i,2}\\
       \vdots  & \ddots  & \vdots & \vdots\\
       0 & \dots & -c_{i,d+1} & c_{i,d+1}
    \end{pmatrix}\in \PGL{d+1}.\]
    Up to conjugation, suppose that $p_0,p_1,\dots,p_{d+1}$ is the standard frame. Then, the monodromy of $P$ is:
    \[M=L_1\dots L_n.\]
\end{lemma}

\begin{proof}
    The proof follows as in $\mathbb{P}^1$. Take the lifts of $p_0,p_1,\dots,p_{d+1},p_{d+2}$:
    \[V_0=\begin{pmatrix}
       1\\
       \vdots\\
       1
    \end{pmatrix},\quad V_1=\begin{pmatrix}
       1\\
       0\\
       \vdots\\
       0
    \end{pmatrix},\quad\dots,\quad V_{d+1}=\begin{pmatrix}
       0\\
       \vdots\\
       0\\
       1
    \end{pmatrix},\quad V_{d+2}=\begin{pmatrix}
       c_{1,1}\\
       \vdots\\
       c_{1,d+1}\\
    \end{pmatrix}.\]
    A straightforward computation gives that
    \[L_1:\begin{cases} 
      p_0 \sim V_0 & \mapsto  c_{1,1} V_1 \sim p_1\\
      p_j \sim V_j & \mapsto  -c_{1,j+1} V_{j+1} \sim p_{j+1}\quad\forall j\in\{1,\dots,d\}\\
      p_{d+1} \sim V_{d+1} & \mapsto  V_{d+2} \sim p_{d+2}
   \end{cases}.\]
   Moreover, because of proposition \ref{prop:cr_inv}, we also have that 
   \[L_1: p_{d+2}\sim V_{d+2}\mapsto V_{d+3}\sim p_{d+3}.\]
    Call $M_i=L_1\dots L_i$. By the exact same recurrence as in \cite[Lemma 3.2]{arnoldCrossratioDynamicsIdeal2022}, we get that 
    \[M_i (p_0,\dots,p_{d+1},c_{i})\mapsto (p_i,\dots,p_{i+d+1},p_{i+d+2}).\]
    Applying it to $M=M_n$ gives us the result.
\end{proof}

Computing the determinants of the $L_i$'s is easy, and yields the following corollary.

\begin{coro}\label{cor:det}
    We have
    \[\mathrm{det}(M)=\prod_{i=1}^n \prod_{j=1}^{d+1} c_{i,j}.\]
\end{coro}

Since a polygonal transformation preserves the monodromy, the coefficients of the characteristic polynomial of $M$ are invariant through the dynamic. They are the usual candidates for Hamiltonians. We compute some of them in section~\ref{sub:preserved}.

\subsection{Scaling symmetry and infinitesimal monodromy}

Now that we have a description of the moduli space in term of cross-ratio coordinates $C:=(c_{i,j})_{1\leq j \leq d+1}^{1\leq i \leq n}$, understanding how a polygonal transformation $T$ acts revolves on finding a formula for $\tilde{C}:=T(C)$. Remember that $T$ is birational, so every $\tilde{c}_{k,l}$ is a rational expression in the $c_{i,j}$'s.\footnote{With this observation, one could say that any polygonal dynamic is a cross-ratio dynamic. But this denomination is already devoted to specific systems.}

Guessing the right equations linking $C$ with $\tilde{C}$ is quite tedious (see lemmas \ref{lemma:flat_cr} and \ref{lemma:stair_cr}). When we get them, we can try to fit a so called \textit{spectral parameter}. This leads to the following definition, more precise than \ref{def:scaling_sym}. It can seem quite convoluted, see propositions \ref{prop:leapfrog_scaling}, \ref{prop:flat_scaling}, \ref{prop:stair_scaling} for examples.

\begin{definition}\label{def:precise_scaling_sym}
    Let $T$ be a polygonal transformation. A \textbf{scaling symmetry} $\rho$ is a non-trivial action of the algebraic group $\mathbb{G}_m$ on $\Tn$,
    \[\rho:\begin{pmatrix}
       \mathbb{G}_m\times\Tn & \to & \Tn\\
       (t,[P]) & \mapsto & [P]_t
    \end{pmatrix},\]
    which commutes with $T$.
\end{definition}

\begin{remark}
    In the cross-ratio coordinates, $\rho$ expresses as a rational map in $C$ and $t$. Hence the action is not defined whenever one of the $c_{i,j}(t)$ is equal to $0$ or $\infty$, which corresponds to a degenerate polygon.
\end{remark}

Finding a scaling symmetry is a little miracle. There isn't, for now, a method to produce them. It is the main tool to find Hamiltonians and prove integrability. It also produces a spectral curve, which is the base ingredient of algebraic integrability.

\begin{remark}\label{remark:P_t}
A scaling symmetry acts on $\Tn$, but it can be lifted to $\tildeTn$. Indeed, if $P$ is an $n$-gon in $\mathbb{P}^d$, we deform its class to get $([P]_t)$ and we lift it by defining
\begin{itemize}
    \item $P_1=P$,
    \item $P_t$ verifies that $[P_t]=[P]_t$, and that $(p_{0}(t),\dots,p_{d+1}(t))$ is constant.   
\end{itemize} 
By the first fundamental theorem of projective geometry, this lift exists and is unique. All the other lifts are given by $(A(t)\cdot P_t)$,  where $(A(t))$ is any algebraic family of projective transformations such that $A(1)= \mathrm{Id}$.
\end{remark}

\begin{remark}
    A closed polygon $P$, when deformed through the scaling symmetry, gives a family of twisted polygons. This is the idea of ``opening up'' a polygon presented in \cite[Fig. 3]{aboudLimitPointPentagram2022}. Thanks to lemma \ref{lemma:monodromy}, we can keep track of the monodromy $M(t)$, whose entries are rational in $C$ and $t$. 
\end{remark}

Since the coefficients of $M(t)$ are rational in $t$, we can formally differentiate it. This gives us our main tool.

\begin{def/prop}\label{prop:infinit_monod}
    Let $\rho$ be a scaling symmetry for a certain polygonal dynamic over $\mathbb{P}^d$.
    
    Let $P$ be a closed polygon and $(P_t)$ be a lift of the deformation by $\rho$. Let $M(t) \in \PGL{d+1}$ be the monodromy of $P_t$ and define the \textbf{infinitesimal monodromy} to be: 
    \[M'(1)\in \mathfrak{pgl}_{d+1}.\]
    It is independent of the choice of lift and invariant under the dynamic. 
\end{def/prop}

\begin{remark}\label{remark:pgl_sum}
    Matrices in $\mathfrak{pgl}_{d+1}$ are defined up to affine transformation. Indeed, for any non-vanishing map $\lambda$, the families of matrices $(M(t))$ and $(\lambda(t)M(t))$ in $\mathrm{GL}_{d+1}$ represent the same projective transformations. Hence, we have:
    \begin{align*}
        (\lambda M)'(1) & =\lambda'(1)\mathrm{Id}+\lambda(1)M'(1).
    \end{align*}
\end{remark}

\begin{proof}
     It is exactly the same proof as in \cite[Prop 3.2]{aboudLimitPointPentagram2022}, which treated the case of $\mathbb{P}^2(\mathbb{C)}$. Let $P_t,~\tilde{P}_t$ be two lifts such that $P_1=\tilde{P}_1$, and $(A(t))$ be an algebraic family of projective transformations such that $\tilde{P}_t=A(t)\cdot P_t$. This forces $A(1)=\mathrm{Id}$ and $\tilde{M}(t)=A(t) M(t)A(t)^{-1}$. By differentiating, we get:
    \[\tilde{M}'(1)= M'(1)+\left[ A'(1),M(1) \right].\]
    Since $P_1$ is a closed polygon, we have $M(1)=\mathrm{Id}$ and so the commutator vanishes. Hence infinitesimal monodromy is well defined.\footnote{As noted in \cite{aboudLimitPointPentagram2022}, the infinitesimal monodromy isn't well defined for twisted polygons, since the only central element of $\PGL{d+1}$ is $\mathrm{Id}$.}
    
    The infinitesimal monodromy is defined from the family of monodromies $(M(t))$ coming from the scaling symmetry. But $T$ preserves the monodromy and commutes with the scaling symmetry, so it also preserves the infinitesimal monodromy.
\end{proof}

The infinitesimal monodromy can be explicitly computed. This is a generalization of Glick's operator, first defined for the pentagram map \cite{glickLimitPointPentagram2020} and then reinterpreted as an infinitesimal monodromy \cite{aboudLimitPointPentagram2022}.

\begin{theorem}\label{thm:infinitesimal_monod}
    Let $T$ be a polygonal transformation in $\mathbb{P}^d$ admitting a scaling symmetry. Let $P=(p_1,\dots,p_n)$ be a closed $n$-gon and $(c_{i,j})$ be its cross-ratio coordinates. Then, its infinitesimal monodromy is
    \[M'(1): p \mapsto \sum_{i=1}^n\sum_{j=0}^{d}\frac{c_{i,j+1}'(1)}{c_{i,j+1}(1)} \frac{\mathrm{det}(p_{i},\dots,p_{i+j-1},p,p_{i+j+1},\dots,p_{i+d})}{\mathrm{det}(p_{i},\dots,p_{i+d})}p_{i+j}.\]
\end{theorem}

\begin{proof}
    We proceed as in the proof of \cite[Prop. 4.1]{aboudLimitPointPentagram2022}, but with our cross-ratio coordinates.
    Let $P$ be a closed polygon. Up to a projective transformation, suppose that $(p_0,\dots,p_d)$ is the standard projective frame. Following remark \ref{remark:P_t}, we deform $P$ into the family $(P_t)$. For any $t\in\mathbb{G}_m$, we lift again the $p_i(t)$'s by choosing non-zero vectors $V_i(t)$'s as in the proof of lemma \ref{lemma:monodromy}. Because of this lemma, we know that for any $t$ the monodromy is 
    \[M(t)=L_1(t)\dots L_n(t).\]
    For any $i\in \{1,\dots,n\}$, define $M_i(t)=L_1(t)\dots L_i(t)$ and $M_0(t)=\mathrm{Id}$. Remark that 
    \[L_{i+1}(t)\dots L_n(t)=M_i(t)^{-1}M(t).\]
    With this, we obtain that 
    \[M'(t)=(\sum_{i=1}^n M_{i-1}(t)L_i'(t)M_i(t)^{-1})M(t).\]
    Since $M(1)=\mathrm{Id}$, this provides 
    \[M'(1)=\sum_{i=1}^n S_i, \qquad \text{where } S_i:=M_{i-1}(1)L_i'(1)M_i(1)^{-1}.\]
    Now we need to understand how each $S_i$ acts. Out of readability, we remove the ``$(1)$''. We have by the proof of lemma $\ref{lemma:monodromy}$ that
    \[M_i:\begin{cases} 
      V_0 & \mapsto  \alpha_{i,0} V_i,\\
      V_j & \mapsto  \alpha_{i,j} V_{i+j} \quad\forall j\in\{1,\dots,d\},\\
      V_{d+1} & \mapsto  \alpha_{i,d+1} V_{i+d+2},
   \end{cases}\]
   where 
   \begin{equation}\label{eq:alphas}
       \alpha_{i,j}=\begin{cases} 
      1 &\text{if } i=0 \text{ or } j>d,\\
      c_{i,1}\alpha_{i-1,1} & \text{if } j=0,\\
      -c_{i,j+1}\alpha_{i-1,j+1} & \text{if } j\in \{1,\dots d\}.
   \end{cases}
   \end{equation}
   We also have 
   \[L_i'=\begin{pmatrix}
       0 & \dots & 0 & c_{i,1}'\\
       -c_{i,2}' & \dots & 0 & c_{i,2}'\\
       \vdots  & \ddots  & \vdots & \vdots\\
       0 & \dots & -c_{i,d+1}' & c_{i,d+1}'
    \end{pmatrix}\]
    so it acts as 
    \[L_i':\begin{cases} 
      V_0 & \mapsto  c_{i,1}' V_1,\\
      V_j & \mapsto  -c_{i,j+1}' V_{j+1} \quad\forall j\in\{1,\dots,d\},\\
      V_{d+1} & \mapsto  c_{i,1}' V_1 + \dots + c_{i,d+1}' V_{d+1}.
   \end{cases}\]
   We assemble these statements to get that $S_i$ acts as:
   \begin{align*}
        V_{i} & \stackrel{M_i^{-1}}{\mapsto} & \frac{1}{\alpha_{i,0}} V_{0} & \stackrel{L_i'}{\mapsto} & \frac{c_{i,1}'}{\alpha_{i,0}} V_{1} & \stackrel{M_{i-1}}{\mapsto} & \frac{c_{i,1}'\alpha_{i-1,1}}{\alpha_{i,0}} V_{i},\\
        V_{i+j} & \mapsto & \frac{1}{\alpha_{i,j}} V_{j} & \mapsto & \frac{-c_{i,j+1}'}{\alpha_{i,j}} V_{j+1} & \mapsto & \frac{-c_{i,j+1}'\alpha_{i-1,j+1}}{\alpha_{i,j}} V_{i+j},
    \end{align*}
    where $j\in\{1,\dots,d\}$. Because of equation \ref{eq:alphas}, we have that the coefficient in front of the output vectors are
    \[\frac{c_{i,j+1}'}{c_{i,j+1}} \quad \text{for all } j\in\{0,\dots,d\}.\]
    Since the polygon $P$ is non degenerate, the family ($V_{i},\dots,V_{i+d}$) forms a basis of $\mathbb{A}^{d+1}$. Hence we can decompose any vector $V$ as
    \[V=\sum_{j=0}^{d} \frac{\mathrm{det}(V_{i},\dots,V_{i+j-1},V,V_{i+j+1},\dots,V_{i+d})}{\mathrm{det}(V_{i},\dots,V_{i+d})}V_{i+j},\]
    and by linearity we obtain 
    \begin{equation*}
        S_i: V\mapsto \sum_{j=0}^{d} \frac{c_{i,j+1}'}{c_{i,j+1}} \frac{\mathrm{det}(V_{i},\dots,V_{i+j-1},V,V_{i+j+1},\dots,V_{i+d})}{\mathrm{det}(V_{i},\dots,V_{i+d})}V_{i+j}.
    \end{equation*}
    We projectivise it back and get that
    \begin{equation}\label{eq:Si}
        S_i: p\mapsto \sum_{j=0}^{d} \frac{c_{i,j+1}'}{c_{i,j+1}} \frac{\mathrm{det}(p_{i},\dots,p_{i+j-1},p,p_{i+j+1},\dots,p_{i+d})}{\mathrm{det}(p_{i},\dots,p_{i+d})}p_{i+j}.
    \end{equation}
    Summing for all $i$'s gives the result.
\end{proof}

\begin{remark}
    As mentioned in remark \ref{remark:pgl_sum}, the infinitesimal monodromy is defined up to an affine transformation. This manifests in theorem \ref{thm:infinitesimal_monod} if we make a different choice of lift for the cross-ratio coordinates, or if we precompose the scaling symmetry $\rho$ with an automorphism of $\mathbb{G}_m$ (which are of the form $t^k$, with $k\in\mathbb{Z}^*$). The new coordinates are given by 
    \[\tilde{c}_{i,j}(t)=\lambda_{i}(t)c_{i,j}(t^k),\]
    where the $\lambda_{i}$'s are non-vanishing functions. Hence we obtain: 
    \begin{align*}
        \frac{\tilde{c}_{i,j}'(1)}{\tilde{c}_{i,j}(1)} & =\frac{\lambda_{i}'(1)c_{i,j}(1)+\lambda_{i}(1)k c_{i,j}'(1)}{\lambda_{i}(1)c_{i,j}(1)}\\
        & = \frac{\lambda_{i}'(1)}{\lambda_{i}(1)}+k\frac{c_{i,j}'(1)}{c_{i,j}(1)}.
    \end{align*}
    Plugging it in equation \ref{eq:Si}, we see that 
    \[\tilde{S}_i=\frac{\lambda_{i}'(1)}{\lambda_{i}(1)}\mathrm{Id}+k S_i,\]
    so by summing over every $i$ we get
    \[\tilde{M}'(1)=\left( \sum_{i=1}^n \frac{\lambda_{i}'(1)}{\lambda_{i}(1)} \right) \mathrm{Id}+k M'(1).\]
\end{remark}

\subsection{Linear scaling symmetry}

We restrain our study to a particular class of scaling symmetries.

\begin{definition}
A scaling symmetry is said to be \textbf{linear} if the normalized cross-ratio coordinates are polynomials of degree at most one in $t$: 
\begin{align*}
    c_{i,j}(t)&=c_{i,j}(1)+c_{i,j}'(1)(t-1) && \forall (i,j)\in \{1,\dots,n\}\times\{1,\dots,d\},\\
    c_{i,d+1}(t)&=1 && \forall i\in \{1,\dots,n\}.
\end{align*}
\end{definition}

So far, all known examples are linear.\footnote{At first glance, this could seem to be wrong for the scaling symmetry of the pentagram map, defined in \cite[Cor. 2.5]{ovsienkoPentagramMapDiscrete2010} as
\[(x_i,y_i)\mapsto(tx_i,t^{-1}y_i).\]
However, the corresponding cross-ratio coordinates should be expressed in terms of $(x_i,y_i^{-1})$, for which the scaling symmetry is linear.} Maybe there is a reason for this. Following the spirit of \cite{arnoldCrossratioDynamicsIdeal2022}, we define for any $I\subset \{1,\dots,n\}\times\{1,\dots,d\}$ the quantities
\[c_I:=\prod_{(i,j)\in I}c_{i,j}(1), \quad c'_I:=\prod_{(i,j)\in I}c'_{i,j}(1)\]

A direct computation gives us the following lemma. 
\begin{lemma}\label{lemma:c_I(t)}
    We have for any $I\subset \{1,\dots,n\}\times\{1,\dots,d\}$ and $t\in\mathbb{G}_m$:
    \[c_I(t)=c_I\sum_{J\subset I}\frac{c_J'}{c_J}(t-1)^{|J|}.\]
\end{lemma}

This gives an simple formula for the determinant.

\begin{lemma}
    Let $T$ be a polygonal transformation in $\mathbb{P}^d$ admitting a linear scaling symmetry. Let $P$ be an $n$-gon (twisted or closed) and $M(t)$ be the deformation of its monodromy from remark \ref{remark:P_t}. For $0\leq k \leq nd$, define:
    \[R_k=\sum_{|I|=k}\frac{c_I'}{c_I}.\]
    They are preserved by the dynamic and projective transformations, and verify:
    \[\mathrm{det}(M(t))=\mathrm{det}(M(1))\sum_{k=0}^{nd}R_k(t-1)^k.\]
\end{lemma}

\begin{proof}
    Corollary \ref{cor:det} tells us that
    \[\mathrm{det}(M(t))=c_{\{1,\dots,n\}\times\{1,\dots,d\}}(t).\]
    Lemma \ref{lemma:c_I(t)} gives us the formula with the $R_k$'s. They are invariant by the dynamic since it preserves the monodromy, and by projective transformations because they are cross-ratios.  
\end{proof}

For closed polygons, we can compute the higher derivatives of the monodromy.

\begin{proposition}\label{prop:higher_derivatives}
    Let $T$ be a polygonal transformation in $\mathbb{P}^d$ admitting a linear scaling symmetry. Let $P$ be a closed $n$-gon.
    Then, for any $k\in\mathbb{N}^*$, the $k$-th infinitesimal monodromy is
    \[M^{(k)}(1)=\sum_{1\leq i_1< \dots < i_k \leq n} S_{i_1}\dots S_{i_k},\]
    where the $S_i$'s are the matrices defined in equation \ref{eq:Si}.
    Again, it is invariant under the dynamic.
\end{proposition}

This provides a lot of preserved quantities, namely all the coefficients of the characteristic polynomial. We'll give a formula for some of them in proposition \ref{prop:Gk}.

\begin{proof}
    The computation is practically the same as the one done for theorem \ref{thm:monod_closed}.
    By the linearity of the scaling symmetry, we know that $L_i''=0$ for any $i$. Hence, we get
    \[M^{(k)}(t)=\sum_{1\leq i_1< \dots < i_k \leq n} L_1(t)\dots L_{i_1-1}(t)L_{i_1}'(t)L_{i_1+1}\dots L_{i_k-1}(t)L_{i_k}'(t)L_{i_k+1}\dots L_{n}(t).\]
    we remark that for any $i<j$:
    \begin{align*}
        L_{i+1}(t)\dots L_{j-1}(t) &= M_i(t)^{-1}M_{j-1}(t).
    \end{align*}
    When we evaluate at $t=1$, we have by the definition of $S_i$:
    \begin{align*}
        M^{(k)}(1) & =\sum_{1\leq i_1< \dots < i_k \leq n} M_{i_1-1}(1)L_{i_1}'(1)M_{i_1}^{-1}(1)\dots M_{i_k-1}(1)L_{i_k}'(1)M_{i_k}^{-1}(1) \\
        & = \sum_{1\leq i_1< \dots < i_k \leq n} S_{i_1}\dots S_{i_k}.
    \end{align*}
    The invariance is again a consequence of the fact that the dynamic preserves the monodromy and commutes with the scaling action.
\end{proof}

\section{Collapsing}\label{section:in P^n}

In any of the dynamics we could simulate with SageMath over $\mathbb{C}$ and $\mathbb{Q}_p$ (unfortunately $\mathbb{C}_p$ is not implemented), we observe a collapse behaviour for the twisted and closed polygons. This is the motivation of the following section.

We'll only consider irreducible polygonal transformations, meaning that they cannot be decomposed in independent smaller systems. A counter-example would be any $2n$-gon $P$ in $\mathbb{P}^2$ composed of intertwining the vertices of two $n$-gons $P_\text{even}$ and $P_{odd}$ with the same monodromy, and applying the pentagram map respectively on them. 

\subsection{For general polygons}

It has been known for a long time that the pentagram map acts trivially on the moduli spaces of closed pentagons and hexagons (see \cite[§2]{schwartzPentagramMap1992}). This phenomenon stops for polygons with at least seven sides. However there exist a family of closed polygons of any number of vertices (satisfying a geometric constrain), whose orbits on the moduli space are periodic (\cite{schwartzPentagramIntegralsPoncelet2015}, \cite{izosimovPentagramMapPoncelet2022}).

This reasoning generalizes. The triviality of the action on the moduli spaces of $(d+3)$- and $(2d+2)$-gons holds for some pentagram maps in $\mathbb{P}^d$ (see \cite{dirdakGALETRANSFORMPENTAGRAM}). Moreover, we prove in proposition~\ref{prop:act_trivial} that any polygonal dynamic in $\mathbb{P}^1$ is trivial on the moduli space of triangles and quadrilaterals, and we construct in subsection \ref{subsub:special} a family of closed $n$-gons (for any $n$) whose dynamic on the moduli space is periodic. More generally, for integrable systems, the invariant manifold of an $n$-gon is a point when $n$ is small enough, hence its dynamic is trivial on the moduli space.

\begin{remark}
    All of the previous examples involve closed polygons. But if the dynamic admits a scaling symmetry, it produces from a polygon $P$ a family of twisted polygons $(P_t)$ whose dynamic is also periodic on the moduli space. This motivates the following theorem.
\end{remark}

\begin{theorem}\label{thm:asympt_periodic}
    Let $F$ be an algebraicly closed field, complete for the metric coming from a non-trivial valuation. Let $T$ be a polygonal transformation taking place in $\mathbb{P}^d(F)$.\\
    Then for almost every $n$-gon $P=((p_1,\dots,p_n),M)$ whose dynamic is periodic on the moduli space, the vertices of $P$ collapse towards a fixed point of $M$.  
\end{theorem}

The proof relies on the following lemmas, which we will prove later. The first one says that iterating a generic element of $\PGL{d+1}(F)$ yields a North-South dynamic, and the second one gives the potential collapsing points.

\begin{lemma}\label{lemma:attract}
    Let $F$ be an algebraicly closed field, complete for the metric coming from a non-trivial valuation
    
    Then a generic element of $\PGL{d+1}(F)$ has one attractive and one repelling fixed point.
\end{lemma}

\begin{lemma}\label{lemma:twist_cv_point}
    Let $P=((p_1,\dots,p_n),M)$ be an $n$-gon. Suppose that under the iteration of $T$, the vertices of P all collapse towards the same value $q$. Then $q$ has to be a fixed point of $M$.
\end{lemma}

Now, let's prove the theorem.

\begin{proof}[Proof of theorem \ref{thm:asympt_periodic}]
    Let $P=((p_1,\dots,p_n),M)$ be an $n$-gon and $N\in \mathbb{N}^*$ be such that $T^N([P])=[P]$. In the initial space, this means that there exists $\tilde{M} \in \PGL{d+1}$ such that $T^N(P)=\tilde{M}\cdot P$. By lemma~\ref{lemma:attract}, we can assume that $\tilde{M}$ almost surely have an attractive and a repelling fixed point, respectively denoted by $q_+$ and $q_-$. Since $T$ commutes with the action of $\PGL{d+1}$, we have:
    \[T^{2N}(P)=T^N(\tilde{M}\cdot P)=\tilde{M}\cdot T^N(P)=\tilde{M}^2\cdot P.\]
    By recurrence (that we can also apply backwards), we get:
    \[l\in\mathbb{Z},~T^{lN}(P)=\tilde{M}^l\cdot P.\]
    We also have that for all $1\leq r<N$ and $l\in\mathbb{{Z}}$:
    \[T^{lN+r}(P)=T^{r}(T^{lN}(P))=T^r(\tilde{M}^l\cdot P)=\tilde{M}^l\cdot T^r(P).\]
    Then
    \[\lim_{l\to \pm \infty} T^{l}(P)=((q_\pm,\dots,q_\pm),M),\]
    and by lemma \ref{lemma:twist_cv_point} the collapse points have to be fixed points of $M$.
\end{proof}

We still have to prove the lemmas.

\begin{proof}[Proof of lemma \ref{lemma:attract}]
    Let $M\in \mathrm{GL}_{d+1}(F)$. The set of diagonalisable matrices of $\mathrm{GL}_{d+1}(F)$ is Zariski dense; indeed, $M$ is diagonalisable if its characteristic polynomial $\chi_M$ have distinct roots. This happens if its discriminant, which is a polynomial in the coefficients $\Delta(a_{1,1},\dots,a_{d+1,d+1})$, doesn't vanish. Since this is not the zero polynomial, we obtain a Zariski-open set. It is not empty since it contains the identity.
    
    So up to conjugation we can assume that:
    \[M=diag(\lambda_1,\dots,\lambda_{d+1}),\]
    with $|\lambda_1|\leq\dots\leq|\lambda_{d+1}|$. Since the valuation is non-trivial, then out of genericity we can assume that the inequalities are in fact strict. In $\PGL{d+1}(F)$, we can take the representative with $\lambda_{d+1}=1$. So:
    \[\lim_{N \to \infty} M^N =E_{d+1,d+1}.\]
    Thus $M$ has an attractive fixed point in $\mathbb{P}^d(F)$, being $[0:\dots:0:1]$ (the origin).
    
    We can do similarly for $M^{-1}$ and take a representative with $\lambda_{1}=1$. This yields:
    \[\lim_{N \to \infty} M^{-N} =E_{1,1},\]
    and the repelling fixed point of $M$ is $[1:0:\dots:0]$ (``the'' point at infinity).
\end{proof}

\begin{proof}[Proof of lemma \ref{lemma:twist_cv_point}]
    Write $P^{(N)}:=T^N(P)$. Since $T$ preserves the monodromy, we have that for all $N\in \mathbb{Z}$:
    \[p_{i+n}^{(N)}=Mp_{i}^{(N)}.\]
    But since $M$ is continuous, by taking the limit we get:
    \[q=Mq. \qedhere \]
\end{proof}

It is known that some polygonal dynamics collapse for many polygons which aren't periodic on the moduli space. The most famous example is the pentagram map on convex polygons \cite[Thm. 3.1]{schwartzPentagramMap1992}. It also happens for other pentagram maps in $\mathbb{P}^2$ \cite[Thm. 1.3]{schwartzFlappingBirdsPentagram2026} and in $\mathbb{P}^d$ (for any $d$) \cite[Thm. 1.4]{yaoGlicksConjecturePoint2014a}. Our computer simulations suggest that the collapsing almost always happens.

\begin{conj}\label{conj:asympt}
    Theorem \ref{thm:asympt_periodic} holds for almost every polygons, even if their dynamic is not periodic on the moduli space.
\end{conj}

If the polygonal dynamic is integrable, then the motion on the moduli space is quasi-periodic. One could try to mimic the case of periodic motion and say that some subsequence of the iterate ``behaves asymptotically the same'' as the iteration of a projective transformation.\footnote{This was already an idea from Schwartz in \cite[§4]{schwartzPentagramMap1992}. He stated in our private communication that no progress have been achieved since 1992.} The meaning of ``behaving asymptotically the same'' is yet unclear. One strategy could be to use generalized schwarzian derivatives (see \cite{ovsienkoProjectiveDifferentialGeometry2004}), but we haven't succeeded.\\
Moreover, computer simulations seem to indicate that the conjecture is true for the projective heat map, which isn't integrable. However, almost every projective classes of pentagons converge towards a fixed point (\cite[Cor. 1.7]{schwartzProjectiveHeatMap2017}) and it is conjecture to be true for any polygon, so they might also behave asymptotically like a projective transformation.

\subsection{For closed polygons}

By lemma \ref{lemma:twist_cv_point}, we know that the only potential collapsing points of a polygon are fixed by its monodromy. However, this is useless for closed polygons (whose monodromy is $\mathrm{Id}$). But if the dynamic admits a scaling symmetry, we can use the infinitesimal monodromy\footnote{See definition \ref{prop:infinit_monod}.} to retrieve the ``true fixed points'' of the identity.
\iffalse
\begin{theorem}\label{thm:collapse_closed}
    Let $T$ be a polygonal transformation taking place in $\mathbb{P}^d(F)$ where $F$ is a complete field, such that the dynamic admits a scaling symmetry. Let $P$ be a closed polygon and $M'$ be its infinitesimal monodromy.
    
    If the vertices of $P$ collapse to a point $q$ under the iteration of $T$, then $q$ has to be an eigenline of $M'$.
\end{theorem}

\begin{proof}
    Suppose that the vertices of $P$ collapse to the point $q$ under the iteration of $T$. Then, since the scaling symmetry conjugates the dynamic, we also have that the vertices of the twisted polygons $P_t=((p_{1}(t),\dots,p_{n}(t)),M(t))$ collapse to some points $q(t)$. Because of lemma~\ref{lemma:twist_cv_point}, they have to be fixed points of $M(t)$. Since the family of matrices $(M(t))$ varies smoothly, we get that the family of collapse points $(q(t))$ also varies smoothly. For $t\sim1$ we have
    \[M(t)=\mathrm{Id}+(t-1)M'+o(t-1),\] so the fixed points of $M(t)$ tend towards the ones of $M'$. Since $q=q(1)$, we get the result.
\end{proof}
\fi
 
\begin{theorem}\label{thm:collapse_closed}
    Let $P$ be a closed polygon satisfying the hypothesis of theorem \ref{thm:asympt_periodic}. Moreover, suppose that the dynamic admits a scaling symmetry. Then, the collapse point $q$ is a fixed point of the infinitesimal monodromy $M'(1)$. 
\end{theorem}

\begin{remark}
    If $M'(1)=0$, the same result holds by replacing it with the first non-vanishing derivative $M^{(k)}(1)$ defined in proposition \ref{prop:higher_derivatives}.
\end{remark}

This result should also be valid if the conjecture \ref{conj:asympt} is true, meaning that the collapse point is always a fixed point of the infinitesimal monodromy.

\begin{proof}
    Because of theorem \ref{thm:asympt_periodic}, we know that the vertices of $P$ converge towards the point $q$. By applying the scaling symmetry, we obtain a family of polygons $P_t=((p_{1}(t),\dots,p_{n}(t)),M(t))$, where the dependence in $t$ is rational (in particular, smooth). Since the scaling symmetry commutes with the dynamic, every $P_t$ is periodic on the moduli space (with the same period $N$). Hence there is a family of projective transformations $(\tilde{M}(t))$ such that 
    \[T^N(P_t)=\tilde{M}(t) \cdot P_t.\]
    The function $\tilde{M}(t)$ is rational in $t$, since it is uniquely determined by the mapping
    \[(p_1(t),\dots,p_{d+2}(t))\mapsto(p_1^{(N)}(t),\dots,p_{d+2}^{(N)}(t)),\]
    and each of these vertices are rational in $t$. Since $\tilde{M}(1)$ generically have an attractive and a repulsive fixed point, it will still be the case for $t$ close enough to $1$. Hence the collapse point $q(t)$ of $P_t$ vary continuously. Moreover, we know by lemma \ref{lemma:twist_cv_point} that the collapse point is a fixed point of $M(t)$. But since 
    \[M(t)=\mathrm{Id}+(t-1)M'+o(t-1),\]
    the fixed points of $M(t)$ tend towards the ones of $M'$ when $t$ goes to $1$. Since $q=q(1)$, we get the result.
\end{proof}

We also get a generalization of \cite[Thm. 1.1]{glickLimitPointPentagram2020}.

\begin{coro}\label{cor:field_ext}
    Let $P$ be a closed $n$-gon in $\mathbb{P}^d$ with vertices  \[(p_i=[p_{i,1}:\dots:p_{i,d+1}])_{i=1,\dots,n}.\]
    Let $q=[q_1:\dots:q_{d+1}]$ be the collapse point predicted by theorem~\ref{thm:collapse_closed}. Then $q_1,\dots,q_{d+1}$ are contained in a field extension of $\mathbb{Q}(p_{1,1},\dots,p_{n,d+1})$, of degree at most $d+1$.
\end{coro}

\begin{proof}
    Since the scaling symmetry is an action of the algebraic group $\mathbb{G}_m$, the coefficients of the infinitesimal monodromy are rational functions in the $p_{i,j}$'s. Because of theorem~\ref{thm:collapse_closed}, the collapse point $q=[q_1:\dots:q_{d+1}]$ is an eigenline of this $(d+1)\times(d+1)$ matrix. To find it, one needs to solve a linear system with $(d+1)$ equations involving the coefficients of the matrix, hence the result.
\end{proof}

%In the case of the pentagram map and its generalizations, we can apply the corollaries. Note that convexity is not involved. This answers the questions raised by Schwartz in \cite[§9.2]{schwartzFlappingBirdsPentagram2026}, at least over $\mathbb{C}$.

Remark that our results don't involve geometric considerations like convexity. They answer the question of Schwartz in \cite[§9.2]{schwartzFlappingBirdsPentagram2026} about the reappearance of Glick's formula some generalized pentagram maps in $\mathbb{P}^2$.

\iffalse
\begin{coro}\label{cor:pentagram}
    Let $F$ be an algebraicly closed complete valued field. For the pentagram map and its generalisations in $\mathbb{P}^2(F)$, the collapse points are fixed points of Glick's operator.
\end{coro}

\begin{proof}
    The generalisation of Glick's lemma to the cousins of the pentagram maps in $\mathbb{P}^2$ is achieved by \cite[Lemma 9.6]{schwartzFlappingBirdsPentagram2026}.
    Aboud and Izosimov proved in \cite{aboudLimitPointPentagram2022} that Glick's operator coincides with the infinitesimal monodromy, and argument in fact holds on any complete field. Combined with theorem~\ref{thm:collapse_closed}, we get the result. 
\end{proof}
\fi

\begin{remark}
    It would be interesting to consider similar statements in the setting of polygonal dynamics on $\mathbb{R}^d$ with the group $\mathrm{Isom}^+(\mathbb{R}^d)$. For instance in $\mathbb{R}^2$, Benoist and Hulin prove in \cite{benoistIterationPliagesQuadrilateres2004} that when $n=4$, the iteration of some element of the folding group is generically conjugated, on the moduli space, to a translation on an elliptic curve. This matches our definition~\ref{def:integr} of integrability. The authors conclude that, on the initial space, the polygons drifts towards a point on the line at infinity, following some ``drifting vector''. When seen on the projective plane, the fixed points of positive isometries are on the line at infinity, which is coherent with theorem \ref{thm:asympt_periodic}.
\end{remark}

\section{Polygonal dynamics in $\mathbb{P}^1$}\label{section:in P^1}

We have more effective statements for polygonal dynamics in $\mathbb{P}^1$. We apply them in the two last sections.

\subsection{The infinitesimal monodromy and the limit point}

Following definition \ref{def:cr_coord}, the cross-ratio coordinates of an $n$-gon $P$ in $\mathbb{P}^1$ is the $n$-periodic sequence
\[c_i=[c_{i,1}:c_{i,2}],\quad i\in\mathbb{Z}.\]
We normalize each $c_{i,2}$ to $1$, and we retrieve the cross-ratios already introduced in \cite{arnoldCrossratioDynamicsIdeal2022} 
\[c_i=[p_i,p_{i+1},p_{i-1},p_{i+2}]=\frac{p_{i}-p_{i-1}}{p_{i}-p_{i+2}}\frac{p_{i+1}-p_{i+2}}{p_{i+1}-p_{i-1}} \in \mathbb{P}^1\setminus\{0,\infty\}.\]

Let's start with a similar result than the one for the pentagram map on the moduli space of small closed polygons (see \cite{schwartzPentagramMap1992}). This might generalize to higher dimensional projective spaces (see \cite{dirdakGALETRANSFORMPENTAGRAM}). 

\begin{proposition}\label{prop:act_trivial}
    Let $T$ be a polygonal dynamic in $\mathbb{P}^1$. Then we have $T=\mathrm{Id}$ on $\mathcal{P}_3$, and $T^2=\mathrm{Id}$ on $\mathcal{P}_4$.
\end{proposition}

\begin{proof}
    Since $\PGL{2}$ is $3$-transitive, then $\mathcal{P}_3$ is a point and $T$ has to act trivially.
    
    Let $P$ be a closed $4$-gon. Because of \cite[Thm. 1]{arnoldCrossratioDynamicsIdeal2022}, the trace of its monodromy in the projective chart $p_0,p_1,p_2=1,\infty,0$ is 
    \[\mathrm{Tr}(M)=1-(c_1+c_2+c_3+c_4)+c_1c_3+c_2c_4.\]
    But we have
    \begin{align*}
        c_2&=[p_2,p_3,p_1,p_0]=(1-c_1),\\
        c_3&=[p_3,p_0,p_2,p_1]=c_1,\\
        c_4&=[p_0,p_1,p_3,p_2]=(1-c_1),
    \end{align*}
    so
    \[\mathrm{Tr}(M)=2c_1(c_1-1).\]
    The first cross-ratio of the iterates $T^N(P)$ always have to be one of the two roots of the polynomial
    \[X^2-X-\frac{\mathrm{Tr}(M)}{2},\]
    which forces $T^2=\mathrm{Id}$ on $\mathcal{P}_4$.
\end{proof}

We have an explicit formula of the infinitesimal monodromy (from theorem \ref{thm:infinitesimal_monod}) and the degree $2$ polynomial equation verified by the limit point (from corollary \ref{cor:field_ext}).

\begin{theorem}\label{thm:monod_closed}
    Let $T$ be a polygonal transformation in $\mathbb{P}^1$ admitting a scaling symmetry. Let $P=(p_1,\dots,p_n)$ be a closed $n$-gon and $(c_{i})$ be its cross-ratio coordinates. Then, its infinitesimal monodromy is 
    \[M'(1)=\frac{n}{2}\mathrm{Id}+\begin{pmatrix}
        J & -K\\ 
        I & -J
    \end{pmatrix},\]
    where
    \[
        I=\sum_{i=1}^n \frac{c_i'(1)}{c_i(1)}\frac{1}{p_i-p_{i+1}},\quad
        J=\frac{1}{2}\sum_{i=1}^n \frac{c_i'(1)}{c_i(1)}\frac{p_i+p_{i+1}}{p_i-p_{i+1}},\quad
        K=\sum_{i=1}^n \frac{c_i'(1)}{c_i(1)}\frac{p_ip_{i+1}}{p_i-p_{i+1}}.
    \]
    Furthermore, its eigenvectors $(X,Y)^T$ satisfy the polynomial equation
    \[\chi_P(X,Y):= IX^2-2JXY+KY^2=0.\]
\end{theorem}

\begin{remark}
    A similar matrix already appeared in \cite[§6.2.1]{arnoldCrossratioDynamicsIdeal2022}. This answers the question asked by Arnold and Arreche at the end of \cite{arnoldSymmediansHyperbolicBarycenters2024}, about the link between the pentagram map and cross-ratio dynamics.
\end{remark}

\begin{remark}
    As we'll see in the proof, the infinitesimal monodromy and the polynomial $\chi_P$ rewrite themselves as:
    \begin{align*}
        M'(1)&=\sum_{i=1}^n \frac{c_i'(1)}{c_i(1)}\frac{1}{p_i-p_{i+1}}\begin{pmatrix}
        p_i & -p_ip_{i+1}\\ 
        1 & -p_{i+1}
    \end{pmatrix},\\
    \chi_P(X,Y)&=\sum_{i=1}^n \frac{c_i'(1)}{c_i(1)}\frac{(X-p_iY)(X-p_{i+1}Y)}{p_i-p_{i+1}}.
    \end{align*}
\end{remark}

\begin{proof}
    Theorem \ref{thm:infinitesimal_monod} tells us the general form of the infinitesimal monodromy. As we said in the beginning of the section, we always normalize $c_{i,2}(t)=1$ (which forces $c_{i,2}'(1)=0$) and identify $c_i(t)$ with $c_{i,1}(t)$ for every $t\in\mathbb{G}_m$. Thus we have
    \[M'(1): p \mapsto \sum_{i=1}^n\frac{c_{i}'(1)}{c_{i}(1)} \frac{\mathrm{det}(p,p_{i+1})}{\mathrm{det}(p_i,p_{i+1})}p_{i}.\]
    We can rewrite it as the matrix
    \[M'(1)=\sum_{i=1}^n \frac{c_i'(1)}{c_i(1)}\frac{1}{p_i-p_{i+1}}\begin{pmatrix}
        p_i & -p_ip_{i+1}\\ 
        1 & -p_{i+1}
    \end{pmatrix},\]
    which we recognize to be 
    \[\frac{n}{2}\mathrm{Id}+\begin{pmatrix}
        J & -K\\ 
        I & -J
    \end{pmatrix}.\]
    Because of remark \ref{remark:pgl_sum}, this represents the same matrix in $\mathfrak{pgl}_2$ than
    \[\begin{pmatrix}
        J & -K\\ 
        I & -J
    \end{pmatrix}.\]
    In both cases, an eigenvector $\begin{pmatrix}
         X\\Y
    \end{pmatrix}$ with eigenvalue $\lambda$ verifies:
    \[\begin{cases}
        JX-KY&=\lambda X\\
        IX-JY&=\lambda Y
    \end{cases}.\]
    By multiplying the first line by $Y$ and the second by $X$, we get by subtracting them:
    \[IX^2-2JXY+KY^2=0.\]
    This concludes the proof.
\end{proof}

\subsection{Preserved quantities}\label{sub:preserved}

We retrieve many preserved quantities. The first ones is a consequence of theorem \ref{thm:monod_closed}.

\begin{coro}
    The quantities $I,J,K$ are preserved by the dynamic, but not by projective transformations. They change under the action of the generators of $\PGL{2}$ in the following way:
    \begin{align*}
        (z \mapsto z-\lambda) \cdot [I,2J,K] & = [I,-(2 I \lambda - 2J) , I \lambda^2 - 2J \lambda + K] \\
        &= [I,- \chi'_{P} (\lambda,1), \chi_{P} (\lambda,1)], \qquad ~\lambda \in F,\\
        (z \mapsto \lambda z) \cdot [I,2J,K] & = [I \lambda^{-1},2J,K \lambda], \qquad\qquad\qquad \lambda \in F^*, \\
        (z \mapsto \frac{-1}{z}) \cdot [I,2J,K] & = [ K , - 2J , I ].
    \end{align*}
    However, the quantity $\Delta=:J^2-IK$ is projective invariant.
\end{coro}

\begin{proof}
    The invariance of $I,J,K$ under the dynamic simply comes from the invariance of the infinitesimal monodromy. The action of projective transformations does not change the cross-ratio coordinates, but only the $p_i$'s. The change of $I,J,K$ under the action of the generators are given by a simple computation.
    
    We have 
    \[\Delta=\begin{vmatrix}
        J & K\\ 
        I & J
    \end{vmatrix},\]
    so it is projective invariant since the determinant is unchanged by conjugation.
\end{proof}

If the scaling symmetry is linear, we have some good candidates for Hamiltonians. Following \cite[§3.1]{arnoldCrossratioDynamicsIdeal2022}, a subset $I\subset\{1,\dots,n\}$ is said to be cyclically sparse if it doesn't contain a pair of consecutive indices (understood modulo $n$). We abbreviate it ``$I \circlearrowleft$ sparse''. 

\begin{proposition}\label{prop:Hk}
    Suppose that the scaling symmetry is linear. Let $P$ be a (twisted or closed) $n$-gon in $\mathbb{P}^1$. 
    For $0\leq k \leq \lfloor n/2 \rfloor$, define
    \[H_k=\sum_{\underset{|I|\geq k}{I \circlearrowleft \text{sparse}}}(-1)^{|I|}c_I\sum_{\underset{|J|= k}{J\subset I}}\frac{c_J'}{c_J}.\]
    These quantities are preserved by the dynamic and are projective invariant. Moreover, we have for every $t\in\mathbb{G}_m$: 
    \[\mathrm{Tr}(M(t))=\sum_{k=0}^{\floor{n/2}}H_k (t-1)^k.\]
\end{proposition}

\begin{proof}
    We know from \cite[Thm. 1]{arnoldCrossratioDynamicsIdeal2022} that for any polygon, the trace of its monodromy is 
    \[\mathrm{Tr}(M(t))=\sum_{k=0}^{\floor{n/2}}(-1)^k\sum_{\underset{|I|= k}{I \circlearrowleft \text{sparse}}}c_I(t).\]
    We obtain the desired formula by applying lemma \ref{lemma:c_I(t)}. The $H_k$'s are invariant under the dynamic because the monodromy is, and under projective transformations because the trace is.
\end{proof}

We also retrieve many projective invariant quantities, which are a generalization of the $G_k$'s defined in \cite[§5.2]{arnoldCrossratioDynamicsIdeal2022}. It would be interesting to have an interpretation of them as alternating semi-perimeters, like in \cite[§5.3]{arnoldCrossratioDynamicsIdeal2022}.

\begin{proposition}\label{prop:Gk}
    Suppose that the scaling symmetry is linear. Let $P$ be a closed $n$-gon in $\mathbb{P}^1$. 
    For $1\leq k \leq \lfloor n/2 \rfloor$, define
    \[G_k=\sum_{i_1<\dots<i_k} \frac{c_{i_1}'(1)}{c_{i_1(1)}}\dots\frac{c_{i_k}'(1)}{c_{i_k}(1)}\frac{(p_{i_1}-p_{i_k+1})(p_{i_2}-p_{i_1+1})\dots(p_{i_k}-p_{i_{k-1}+1})}{(p_{i_1}-p_{i_1+1})(p_{i_2}-p_{i_2+1})\dots(p_{i_k}-p_{i_{k}+1})}\]
    and $G_0=2$. These quantities are preserved by the dynamic and are projective invariant. Moreover, we have for every $t\in\mathbb{G}_m$:
    \[\mathrm{Tr}(M(t))=\sum_{k=0}^{\floor{n/2}} \frac{G_k}{k!}(t-1)^k.\]
\end{proposition}

Combining the two properties tells us that for every $0\leq k \leq \floor{n/2}$, the quanties are related by 
\[H_k=\frac{G_k}{k!}.\]

\begin{proof}
    By proposition \ref{prop:higher_derivatives}, we have that for any $k\in\mathbb{N}$:
    \begin{align*}
        M^{(k)}(1) & =\sum_{1\leq i_1< \dots < i_k \leq n} S_{i_1}\dots S_{i_k}\\
        & = \sum_{1\leq i_1< \dots < i_k \leq n}\prod_{l=1}^k \frac{c_{i_l}'(1)}{c_{i_l}(1)}\frac{1}{p_{i_l}-p_{i_l+1}}\begin{pmatrix}
        p_{i_l} & -p_{i_l}p_{i_l+1}\\ 
        1 & -p_{i_l+1}
    \end{pmatrix}.
    \end{align*}
    But by the same computation as in the proof of \cite[Thm. 10]{arnoldCrossratioDynamicsIdeal2022}, we obtain that $M^{(k)}(1)=0$ if $k>\floor{n/2}$, and otherwise
    \[\mathrm{Tr}(M^{(k)}(1))=G_k.\]
    We conclude with Taylor's theorem. Their invariance is inherited from the one of the $H_k$'s.
\end{proof}

Let's finish this section with some questions.

\begin{openquestion}
    Is there a automatic way to produce an invariant Poisson bracket for which the $H_k$'s are the Hamiltonians, which could prove Liouville integrability as in \cite{arnoldCrossratioDynamicsIdeal2022} ?   
\end{openquestion}

\begin{openquestion}
    Define the following 1-form and (usually presymplectic) 2-form on the space of closed $n$-gons $\tildePn$:
    \[\lambda=\frac{1}{2}\sum_{i=1}^n \frac{c_i'(1)}{c_i(1)}\frac{dp_i + dp_{i+1}}{p_{i+1}-p_i}, \quad \Omega=d\lambda=\sum_{i=1}^n \frac{c_i'(1)}{c_i(1)} \frac{dp_i \wedge dp_{i+1}}{(p_i-p_{i+1})^2}.\]
    Consider the generators of the infinitesimal action of $\PGL{2}$
    \begin{align*}
        u&=\sum_{1\leq k\leq n}\frac{\partial}{\partial p_{k}}, & v&=\sum_{1\leq k\leq n}p_{k}\frac{\partial}{\partial p_{k}}, & w=&\sum_{1\leq k\leq n}p_{k}^2\frac{\partial}{\partial p_{k}}.
    \end{align*}
    A computation similar to the one done in \cite[§6.2]{arnoldCrossratioDynamicsIdeal2022} shows that 
    \begin{align*}  
        i_u \Omega&=dI, & i_v \Omega&=dJ, & i_w \Omega&=dK +\sum_{i=1}^n \frac{1}{2}(\frac{c_i'(1)}{c_i(1)}-\frac{c_{i-1}'(1)}{c_{i-1}(1)})dp_i.
    \end{align*}
    If all the $\frac{c_i'(1)}{c_i(1)}$ are equal,\footnote{For the systems considered here, it only happens for the flat cross-ratio dynamics with constant discrete curvature.} then the last sum vanishes and this tells us that $I,J,K$ are the hamiltonians of respectively $u,v,w$.
    
    Long specific computations show that in each system, $\Omega$ is invariant through the dynamic, meaning $T^*\Omega=\Omega$ (see \cite[Thm. 15]{arnoldCrossratioDynamicsIdeal2022} and \cite[Cor. 5.15]{gekhtmanIntegrableClusterDynamics2016}). This fact seems to be true for any polygonal dynamic in $\mathbb{P}^1$. An idea to prove it would be to rephrase the dynamic as a cluster transformation and generalize the interpretation of $\lambda$-length made in \cite[§5.3]{arnoldCrossratioDynamicsIdeal2022} with our cross-ratio coordinates.\footnote{This is coherent with the fact that 
    \[\frac{c_i'(1)}{c_i(1)}=\left.\frac{d\mathrm{ln}(c_i(t))}{dt}\right|_{t=1}.\]} Then, $\Omega$ could be the Weil–Petersson Kähler two-form presented in
    \cite[§2.3]{pennerDecoratedTeichmullerTheory2012}.
\end{openquestion}

\section{Two examples of dynamics in $\mathbb{P}^1$}\label{section:examples}
We apply the results of the previous section to three kinds of polygonal dynamics: the leapfrog map, the ``flat'' cross-ratio and the newly defined ``stair-case'' cross-ratio dynamics. The two last ones are linked to the theory of discrete holomorphic functions (see \cite[§6]{affolterSchwarzianOctahedronRecurrence2024}) and are endowed with a \textbf{discrete curvature}, which is an $n$-periodic sequence $\mu=(\mu_i)_{i\in \mathbb{Z}}$, where each $\mu_i\in \mathbb{G}_m$ should be thought to live on the edge linking $p_i$ and $p_{i+1}$. They are left unchanged by the action of $\PGL{2}$, but can be modified by the scaling symmetry.
\subsection{The leapfrog map}

The leapfrog map is defined in \cite[§5.2]{gekhtmanIntegrableClusterDynamics2016} as an analogue of the pentagram map in $\mathbb{P}^1$, and the authors notice that it seems to be the hyperbolic counterpart of the polygon recutting defined in \cite{adlerRecuttingsPolygons1993} (which would be the euclidean case). For a matter of internal coherence, the notations of \cite{gekhtmanIntegrableClusterDynamics2016} are changed to fit the ones from here.

The dynamic involves two $n$-gons
\[S^-=(s_i^-)_{i\in\mathbb{Z}},\quad S=(s_i)_{i\in\mathbb{Z}},\]
with the same monodromy, which we intertwine this way
\[\dots,s_i^-,s_i,s_{i+1}^-,s_{i+1},\dots\]
to form a $2n$-gon denoted by $P$. The leapfrog map is defined by 
\[\Phi:(S^-,S)\mapsto(S,S^+),\]
where $S^+$ is given by applying local ``leapfrog'' moves. Such a move maps $s_i^-$ to $s_i^+$ via the only projective transformation that fixes $s_i$ and interchanges $s_{i-1}$ with $s_{i+1}$. Note that such moves are involutions, and that we could study the group generated by them, like we will do afterwards for the staircase cross-ratio dynamics. However, it is not the point here.

Because of the indexation, it is more convenient work with
\[o_i^-=[s_i^-,s_{i},s_{i-1},s_{i+1}^-], \quad o_i=[s_i,s_{i+1}^-,s_{i}^-,s_{i+1}],\]
instead of the usual $c_i$. This gives a fairly simple scaling symmetry.

\begin{proposition}\label{prop:leapfrog_scaling}
    There is a scaling symmetry acting as
    \[o_i^-(t)=to_i^-, \qquad o_i(t)=o_i. \]
    It is well defined for any $t\in\mathbb{G}_m$.
\end{proposition}

\begin{proof}
    First of all, it is direct to check that this is indeed a group action from $\mathbb{G}_m$.
    
    We express the system of coordinates $(o^-,o)$ with the help of two other systems, and then get the result. According to \cite[Prop. 5.12]{gekhtmanIntegrableClusterDynamics2016}, we have the coordinate systems $(x,y)$ and $(p,r)$, defined by:
    \begin{align*}
        x_i&=[s_{i+1},s_{i+1}^-,s_{i+2}^-,s_{i}^-],&y_i&=\frac{(s_{i+1}^--s_{i+1})(s_{i+2}^--s_{i+2})(s_{i}^--s_{i+1}^-)}{(s_{i+1}^--s_{i+2})(s_{i}^--s_{i+1})(s_{i+1}^--s_{i+2}^-)},\\
        p_i&=[s_{i+1}^-,s_{i+2}^-,s_{i+1},s_{i+2}],& r_i&=[s_i,s_{i+1},s_{i+1}^-,s_{i+2}^-].
    \end{align*}
    They are linked by the relations\footnote{The original paper forgot the minus signs, but it doesn't impact their reasoning.}
    \begin{align*}
        p_i&=-\frac{y_i}{x_i}, & r_i&=-\frac{x_{i-1}x_i}{x_{i-1}(1-x_i)+y_{i-1}}.
    \end{align*}

    Now, because of the permutation of the variables of the cross-ratio and the previous relation, we see that:
    \begin{align*}
        o_i^-&=\frac{r_{i-1}}{r_{i-1}-1}
        =\frac{x_{i-2}x_{i-1}}{x_{i-2}+y_{i-2}},&
        o_i&=\frac{p_{i-1}}{p_{i-1}-1}
        =\frac{y_{i-1}}{x_{i-1}+y_{i-1}}.
    \end{align*}
    But according to \cite[Rmk. 3.5]{gekhtmanIntegrableClusterDynamics2016}, there is a scaling symmetry (actually working for all of their generalised pentagram maps) defined for any $t\in\mathbb{G}_m$ by
    \begin{align*}
        x_i(t)&=tx_i,&  y_i(t)&=ty_i. 
    \end{align*}
    Plugging it back to our formula linking $(o^-,o)$ to $(x,y)$, we get the result.
\end{proof}

Now we derive the infinitesimal monodromy and potential collapse points for the leapfrog map.

\begin{proposition}
    The infinitesimal monodromy for the leapfrog map is
    \[M'(1)=\sum_{i=1}^n \frac{1}{s_i^--s_{i}}\begin{pmatrix}
        s_i^- & -s_i^-s_{i}\\ 
        1 & -s_{i}
    \end{pmatrix},\]
    and its eigenvectors $(X,Y)^T$ are roots of the polynomial 
    \[\chi_P(X,Y)=\sum_{i=1}^n \frac{(X-s_i^-Y)(X-s_{i}Y)}{s_i^--s_{i}}.\]
\end{proposition}

\begin{proof}
    We have from proposition \ref{prop:leapfrog_scaling} that 
    \[
        \frac{o_i^{-}~'(1)}{o_i^-}=1, \qquad \frac{o_i'(1)}{o_i}=0,
    \]
    and by remembering that they correspond to the cross-ratio coordinates $c_i$, we get the result by applying theorem~\ref{thm:monod_closed}.
\end{proof}

\subsection{Flat cross-ratio dynamics}\label{sub:flat_cr}

The flat cross-ratio dynamics arise from the continuation of discrete holomorphic functions, for which the geometry of the initial condition (meaning $P$) is flat. See figure \ref{fig:flat_crossratio} and compare with the staircase cross-ratio dynamics depicted on figure \ref{fig:stair_crossratio}.

It was studied in \cite{arnoldCrossratioDynamicsIdeal2022} and \cite{hetrich-jerominPeriodicDiscreteConformal2001} for constant discrete curvature $\alpha\in\mathbb{C}\setminus\{0,1\}$, where they respectively prove Hamiltonian and algebraic integrability. It was generalized for varying $(\alpha_i)\in(\mathbb{C}\setminus\{0,1\})^n$ in \cite{affolterIntegrableDynamicsProjective2023}.

\begin{definition}
    Let $n\geq3$ and $P=((p_1,\dots,p_n),M)$ be an $n$-gon in $\mathbb{P}^1$. Enrich it with a discrete curvature $\alpha=(\alpha_1,\dots,\alpha_n)$ such that $\alpha_i\in\mathbb{G}_m\setminus\{1\}$ for each $i\in \mathbb{Z}$.
    
    We say that two $n$-gons $P$ and $Q$ are $\alpha$-related if they share the same monodromy, the same discrete curvature, and verify    \[\alpha_i=[p_i,p_{i+1},q_i,q_{i+1}] \quad \forall i \in \mathbb{Z}.\]
    In a more compact way, we write $P \overset{\alpha}{\sim} Q$.
    
    From \cite[Cor. 2.7]{arnoldCrossratioDynamicsIdeal2022} we get that for an $n$-gon $P$, there are generically only two $n$-gon $Q^+,Q^-$ that are $\alpha$-related to it. This gives an evolution operator in the futur/past, called the \textbf{flat cross-ratio dynamics}.
\end{definition}

\begin{figure}
    \centering
    \begin{tikzpicture}
        \draw[gray!30] (0,0) grid (7,2);
        \foreach \i in {0,1,...,7}{
            \foreach \j in {0,1,...,2}{
                %\pgfmathsetmacro\h{int(\i+\j)};
                \filldraw[black] (\i,\j) circle(.8pt);
            }
        }
        \node[left] at (0,1) (0) {$\dots$};
        \node[right] at (9,1) (0) {$\dots$};
        
        \node at (-1,1) (0) {$P$};
        \node at (-1,2) (0) {$Q^+$};
        \node at (-1,0) (0) {$Q^-$};

        \foreach \i in {1,2,...,6}{
            \node[right] at (\i,1.25) (0) {$p_\i$};
            \node[right] at (\i,2.25) (0) {$q_\i^+$};
            \node[right] at (\i,0.25) (0) {$q_\i^-$};
            \node[below, blue] at (\i+0.5,-0.5)  {$\alpha_\i$};
            \draw[dotted, blue] (\i+0.5,-0.5) -- (\i+0.5,2.5);
        }
        
        \node[right] at (0,1.25) (0) {$p_0$};
        \node[right] at (0,2.25) (0) {$q_0^+$};
        \node[right] at (0,0.25) (0) {$q_0^-$};
        \node[below, blue] at (0+0.5,-0.5)  {$\alpha_6$};
        \draw[dotted, blue] (0+0.5,-0.5) -- (0.5,2.5);
    
        \node[right] at (7,1.25) (0) {$p_7=M(p_1)$};
        \node[right] at (7,2.25) (0) {$q_7^+=M(q_1^+)$};
        \node[right] at (7,0.25) (0) {$q_7^-=M(q_1^-)$};
        
        \draw (0,0) -- (7,0);
        \draw (0,1) -- (7,1);
        \draw (0,2) -- (7,2);
        \end{tikzpicture} 
    \caption{Flat cross-ratio dynamics seen in the context of discrete holomorphic functions. The $6$-gon $P=((p_1,\dots,p_6),M)$ is $\alpha$-related in the past to $Q^-$ and in the future to $Q^+$.}
    \label{fig:flat_crossratio}
\end{figure}

Next lemma is an adaptation of \cite[Lemma 4.4]{arnoldCrossratioDynamicsIdeal2022} for non-constant $\alpha_i$'s.

\begin{lemma}\label{lemma:flat_cr}
    Let $P$ and $Q$ be two $n$-gons such that $P \overset{\alpha}{\sim} Q$. For all $i\in \mathbb{Z}$, set
    \begin{align*}
        c_i&=[p_i, p_{i+1}, p_{i-1}, p_{i+2}],& x_i&=[p_i, p_{i+1}, p_{i-1}, q_{i}],\\ 
        d_i&=[q_i, q_{i+1}, q_{i-1}, q_{i+2}],& y_i&=[q_i, q_{i+1}, q_{i-1}, p_{i}].
    \end{align*}
    Then we have 
    \begin{equation}\label{eq:flat_cr}
        \begin{split}
        c_i = \alpha_ix_i(1 - x_{i+1}),\quad d_i = \alpha_i y_i(1 - y_{i+1}),\\\quad 1=x_i + y_i+x_iy_i\left(\frac{\alpha_i-1}{\alpha_{i-1}-1}-1\right).
        \end{split}
    \end{equation}
   
    Conversely, for any $(x_i) \in (\mathbb{G}_m \setminus \{ 1\})^n$  and any $(\alpha_i) \in (\mathbb{G}_m \setminus \{ 1\})^n$, the $(c_i)$ and $(d_i)$ given by the previous formulas define a pair $P \overset{\alpha}{\sim} Q$, unique up to a simultaneous projective transformation.
\end{lemma}

\begin{proof}
    Again, it is similar to \cite{arnoldCrossratioDynamicsIdeal2022}. By using ``Chasles' relation of the cross-ratio'', we have
        \begin{align*}
            c_i&=[p_i, p_{i+1}, p_{i-1}, q_i][p_i, p_{i+1}, q_i, q_{i+1}][p_i, p_{i+1}, q_{i+1}, p_{i+2}]\\
            &= x_i \alpha_i(1 - x_{i+1}),
        \end{align*}
    and similarly we get the same one for $d_i$ by exchanging the $p_i$'s with $q_i$'s, and $x_i$'s with $y_i$'s. We also have
    \begin{align*}
        [p_i, q_i, p_{i-1}, p_{i+1}] &= [p_i, q_i, p_{i-1}, q_{i-1}][p_i, q_i, q_{i-1}, q_{i+1}][p_i, q_i, q_{i+1}, p_{i+1}]\\
        &=(1-\alpha_{i-1})[p_i, q_i, q_{i-1}, q_{i+1}]\frac{1}{1-\alpha_{i}}
    \end{align*}
    so by referring to the permutations of the variables of the cross-ratio, we get
    \[\frac{x_i}{1-x_i}(1-\alpha_i)=\frac{1-y_i}{y_i}(1-\alpha_{i-1}).\]
    This rewrites as
    \[1=x_i + y_i+x_iy_i\left(\frac{\alpha_i-1}{\alpha_{i-1}-1}-1\right).\]
    
    Conversely, the data of $x$ and $\alpha$ allows to reconstruct $y$, $c$ and $d$. We extend them to get $n$-periodic sequences. From $c$, we can choose representative $P=(p_i)_{i\in\mathbb{Z}}$ such $c_i=[p_i,p_{i+1},p_{i-1},p_{i+2}]$. Because of lemma~\ref{lemma:monodromy}, we know its monodromy $M_P$. Similarly we construct a polygon $Q=(q_i)_{i\in\mathbb{Z}}$, with monodromy $M_Q$. We observe that
    \[[p_i, p_{i+1}, p_{i-1}, q_{i}]=x_i=x_{i+n}=[M_P(p_i), M_P(p_{i+1}), M_P(p_{i-1}), M_Q(q_{i})],\]
    hence we have $M_P=M_Q$. Because of the computation that gave us
    \[c_i = \alpha_ix_i(1 - x_{i+1}),\]
    we recognize that $[p_i,p_{i+1},q_i,q_{i+1}]=\alpha_i$, so we have indeed $P \overset{\alpha}{\sim} Q$.
\end{proof}

This crucial lemma allows us to define a scaling symmetry. This was already observed in \cite[Lemma 4.4]{arnoldCrossratioDynamicsIdeal2022} for constant $\alpha$, although the scaling was much simpler and the concept of scaling symmetry was not named.

\begin{proposition}\label{prop:flat_scaling}
    We have a scaling symmetry $\rho$ that acts both on the $c_i$'s and the $\alpha_i$'s, in the following way:
    \[\alpha_i(t)=1+\kappa_i(t\alpha_1-1), \quad c_i(t)=c_i\frac{\alpha_i(t)}{\alpha_i},\qquad \text{where }\kappa_i=\frac{\alpha_i-1}{\alpha_1-1}.\]
    This is well defined for
    $t\neq\frac{1}{\alpha_1},\frac{\kappa_i-1}{\kappa_i \alpha_1}$.
\end{proposition}

In any case, this scaling symmetry is well defined around $t=1$. Indeed, $1/\alpha_1 \neq 1 $ and $(\kappa_i-1)/\kappa_i \alpha_1=1$ would imply $\alpha_i=0$, which is not the case.

\begin{proof}
    Let's check first that $\rho$ is indeed a group action. First, we have
    \[\alpha_i(1)=1+\frac{\alpha_i-1}{\alpha_1-1}(\alpha_1-1)=\alpha_i,\]
    and so $c_i(1)=c_i$. Furthermore, we note that 
    \[\kappa_i(t)=\frac{\alpha_i(t)-1}{\alpha_1(t)-1}=\frac{\kappa_i(t\alpha_1-1)}{(t\alpha_1-1)}=\kappa_i,\]
    so $\kappa_i(t)$ is constant for each $i$ ($\rho$ was in fact designed from that property). Then we have
    \begin{align*}
        \alpha_i(t_1)(t_2)&=1+\kappa_i(t_2\alpha_1(t_1)-1)\\
        &=1+\kappa_i(t_2t_1\alpha_1-1)\\
        &=\alpha_i(t_1t_2),
    \end{align*}
    and from this we obtain
    \begin{align*}
    c_i(t_1)(t_2)&=c_i(t_1)\frac{\alpha_i(t_1)(t_2)}{\alpha_i(t_1)}\\
    &=c_i\frac{\alpha_i(t_1)}{\alpha_i}\frac{\alpha_i(t_1t_2)}{\alpha_i(t_1)}\\
    &=c_i(t_1t_2).
    \end{align*}
    So $\rho$ is indeed a group action. It is well defined if $\alpha_i(t)\neq0,1$ for each $i$, which corresponds to $t\neq\frac{1}{\alpha_1},\frac{\kappa_i-1}{\kappa_i \alpha_1}$.

    Let's now see why it commutes with the evolution operator. Take $c,d,x,y$ as in lemma~\ref{lemma:flat_cr}. They are linked by equation \ref{eq:flat_cr}. Note that $\frac{\alpha_i-1}{\alpha_{i-1}-1}$ is equal to $\frac{\kappa_i}{\kappa_{i-1}}$, so it stays constant when applying $\rho$. In fact, applying $\rho$ modifies $c,d,\alpha$, but leaves $x$ and $y$ unchanged. Equation \ref{eq:flat_cr} still holds for $c_i(t),d_i(t),\alpha_i(t),x_i,y_i$ and so by lemma~\ref{lemma:flat_cr} we have $P(t) \overset{\alpha(t)}{\sim} Q(t)$. Hence $\rho$ is indeed a scaling symmetry.
\end{proof}

\begin{proposition}\label{prop:flat_infinitesimal}
    The infinitesimal monodromy for the flat cross-ratio dynamics is
    \[M'(1)=\sum_{i=1}^n \frac{(\alpha_i-1)\alpha_1}{(\alpha_1-1)\alpha_i}\frac{1}{p_i-p_{i+1}}\begin{pmatrix}
        p_i & -p_ip_{i+1}\\ 
        1 & -p_{i+1}
    \end{pmatrix},\]
    and its eigenvectors $(X,Y)^T$ are roots of the polynomial 
    \[\chi_P(X,Y)=\sum_{i=1}^n \frac{(\alpha_i-1)}{\alpha_i}\frac{(X-p_iY)(X-p_{i+1}Y)}{p_i-p_{i+1}}.\]
\end{proposition}

\begin{proof}
    From proposition \ref{prop:flat_scaling} we see that 
    \[\frac{c_i'(1)}{c_i(1)}=\frac{(\alpha_i-1)\alpha_1}{(\alpha_1-1)\alpha_i},\]
    so by theorem~\ref{thm:monod_closed} we get the formula for the infinitesimal monodromy. Still from the same theorem, we obtain that an eigenvector $(X,Y)^T$ is a root of the polynomial $\chi_P(X,Y)$. Up to multiplying it by $\frac{\alpha_1-1}{\alpha_1}$, which we do only for cosmetic purpose, we get desired formula for the polynomial.
\end{proof}

\section{A group dynamic in $\mathbb{P}^1$: staircase cross-ratio dynamics}\label{section:staircase}

\begin{figure}
    \centering
    \begin{tikzpicture}
        \draw[gray!30] (0,0) grid (6,1);
        \foreach \i in {0,1,...,6}{
            \foreach \j in {0,1,...,1}{
                %\pgfmathsetmacro\h{int(\i+\j)};
                \filldraw[black] (\i,\j) circle(.8pt);
            }
        }
        \node[left] at (0,1) (0) {$\dots$};
        \node[right] at (8,0) (0) {$\dots$};
        
        \node at (-1,1) (0) {$P$};

        \foreach \i in {1}{
            \node[right] at (\i,1.25) (0) {$p_\i$};
            \node[below, blue] at (\i+0.5,-0.5)  {$\mu_\i$};
            \draw[dotted, blue] (\i+0.5,-0.5) -- (\i+0.5,1.5);
        }
        \node[right] at (2,1.25) (0) {$p_2$};
        \node[left,blue] at (-0.5,0.5) (0) {$\mu_2$};
        \draw[dotted, blue] (-0.5,0.5) -- (6.5,0.5);
        
        \node[right] at (3,1.25) (0) {$\tilde{p}_3$};
        %\draw[line width=1pt, dotted,->] (2,0) -- (3,1);
        \draw[line width=1pt, dotted] (2,1) -- (3,1);
        \draw[line width=1pt, dotted] (3,0) -- (3,1);

        \foreach \i in {3,4,...,6}{
            \node[right] at (\i-1,-0.25) (0) {$p_\i$};
            \node[below, blue] at (\i-0.5,-0.5)  {$\mu_\i$};
            \draw[dotted, blue] (\i-0.5,-0.5) -- (\i-0.5,1.5);
        }
        
        \node[right] at (0,1.25) (0) {$p_0$};
        
        \node[right] at (6,-0.25) (0) {$p_7=M(p_1)$};
        
        \node[below, blue] at (0+0.5,-0.5)  {$\mu_6$};
        \draw[dotted, blue] (0+0.5,-0.5) -- (0.5,1.5);
    
        \draw (2,0) -- (6,0);
        \draw (2,0) -- (2,1);
        \draw (0,1) -- (2,1);
        \end{tikzpicture} 
    \caption{Staircase cross-ratio dynamics seen in the context of discrete holomorphic functions. Applying the flip $\phi_3$ on the $6$-gon $P=((p_1,\dots,p_6),M)$ changes $p_3$ to $\tilde{p}_3$ and interchanges $\mu_2$ with $\mu_3$.}
    \label{fig:stair_crossratio}
\end{figure}

We define a new polygonal dynamic in $\mathbb{P}^1$ given by a group of polygonal transformations called flips.\footnote{It most certainly admits an interpretation as cluster transformations.} The staircase cross-ratio dynamics arise from the continuation of discrete holomorphic functions, for which the geometry of the initial condition (meaning $P$) is staircase-like. See figure \ref{fig:stair_crossratio} and compare with the flat cross-ratio dynamics depicted on figure \ref{fig:flat_crossratio}.

However, this dynamic is more general than the problem of continuation of discrete holomorphic functions, for which some flips are not allowed by the geometry.\footnote{On the example presented on figure \ref{fig:stair_crossratio}, the flip $\phi_5$ is not allowed since $p_5$ is not on a ``corner'' like $p_3$.} The more abstract way to consider the dynamic is to forget the underlying geometry and to see the group as a braid group, as on figure \ref{fig:braid} and theorem \ref{thm:affine_sym_grp}.

\subsection{Definition and results}

Before defining the group, let's start with a few preliminaries.

\begin{definition}
    Let $\alpha, \beta \in \mathbb{G}_m$. Two couples of distinct points $(a,c), (b,d)$ of $\mathbb{P}^1$ form a \textbf{$\bm{\alpha/\beta}$-quadrangle} when $[a,c,b,d]=\alpha/\beta$. The point $c$ is called the \textbf{$\bm{\alpha/\beta}$-conjugate} of $a$ with respect to $(b,d)$, and we write \[c=h^{\alpha/\beta}_{b,d}(a).\]
    When $\alpha/\beta=-1$, the quadrangle said to be harmonic.
\end{definition}

\begin{proposition}\label{prop:harm}The following statements are true.
\begin{enumerate}
        \item Projective transformations send $\alpha/\beta$-quadrangles to $\alpha/\beta$-quadrangles. Moreover, we have
        \[ g\circ h^{\alpha/\beta}_{b,d}= h^{\alpha/\beta}_{g(b),g(d)}\circ g, \quad \text{for all }g \in \PGL{2}.\]
        \item The $\alpha/\beta$-conjugation with respect to $(b,d)$ is a projective transformation 
        \[h^{\alpha/\beta}_{b,d}(z)=\frac{(\beta d-\alpha b)z-(\beta-\alpha)bd}{(\beta-\alpha)z-(\beta b-\alpha d)},\]
        which fixes $b$ and $d$, and has characteristic constant $\alpha/\beta$. Its inverse is $h^{\alpha/\beta}_{d,b}$.
        \item If $(b,d)$ is fixed, then for $\alpha,\beta,\gamma,\delta \in \mathbb{G}_m$ we have 
        \[h^{\alpha/\beta}_{b,d}\circ h^{\gamma/\delta}_{b,d}=h^{\alpha\gamma/\beta\delta}_{b,d}.\]
        Hence $h^1_{b,d}=\mathrm{Id}$ and $(h^{\alpha/\beta}_{b,d})^{-1}=h^{\beta/\alpha}_{b,d}$.
    \end{enumerate}
\end{proposition}

\begin{proof}Let's check each statement.
   \begin{enumerate}%[label=(\alph*)]
        \item Let $(a,c), (b,d)$ be an $\alpha/\beta$-quadrangle and $g\in\PGL{2}$. The cross-ratio is $\PGL{2}$-invariant, so 
        \[[g(a),g(h^{\alpha/\beta}_{b,d}(a)),g(b),g(d)]=\frac{\alpha}{\beta}=[g(a),h^{\alpha/\beta}_{g(b),g(d)}( g(a)),g(b),g(d)],\]
        which gives the result.
        \item The formula of $h^{\alpha/\beta}_{b,d}$ is the result of direct calculation. To see the characteristic constant, it suffices to conjugate by a projective transformation that sends $(b,d)$ to $(\infty,0)$, use \textit{(a)} and by the formula we see that $h^{\alpha/\beta}_{\infty,0}:z\mapsto \frac{\alpha}{\beta}z$.\\
        Considering $(c,a),(d,b)$ also gives a $\alpha/\beta$-quadrangle, because this permutation doesn't modify the cross-ratio. So we get \[a=h^{\alpha/\beta}_{d,b}(c)=h^{\alpha/\beta}_{d,b}\circ h^{\alpha/\beta}_{b,d}(a),\]
        hence $h^{\alpha/\beta}_{d,b}=(h^{\alpha/\beta}_{b,d})^{-1}$.
        \item Let $g \in \PGL{2}$ which sends $(b,d)$ to $(\infty,0)$. Using \textit{(a)} and \textit{(b)}, we obtain
        \begin{align*}
            h^{\alpha/\beta}_{b,d}\circ h^{\gamma/\delta}_{b,d}&=g^{-1}\circ h^{\alpha/\beta}_{\infty,0} \circ g \circ g^{-1}\circ h^{\gamma/\delta}_{\infty,0} \circ g \\
            &= g^{-1}\circ h^{\alpha\gamma/\beta\delta}_{\infty,0} \circ g \\
            &= h^{\alpha\gamma/\beta\delta}_{b,d}. \qedhere
        \end{align*}     
    \end{enumerate}
\end{proof}
We can now define the dynamic. See figure~\ref{fig:braid} for a visualisation of it.
\begin{definition}
Let $n\geq3$ and $P=((p_1,\dots,p_n),M)$ be an $n$-gon in $\mathbb{P}^1$. Enrich it with an $n$-periodic discrete curvature $\mu=(\mu_i)_{i\in\mathbb{Z}}$ such that $\mu_i\in\mathbb{G}_m$ for each $i\in \mathbb{Z}$.

Let $j\in\Zn$. The \textbf{flip at index} $\bm{j}$ is the map
\[\phi_j:(P,\mu)\mapsto(\tilde{P},\tilde{\mu})\]
such that:
\[\begin{cases}
    \tilde{p}_i=h^{\mu_i/\mu_{i-1}}_{p_{i+1},p_{i-1}}(p_{i})&\quad\forall i\in n\mathbb{Z}+j,\\
    \tilde{p}_i=p_i&\quad\forall i\notin n\mathbb{Z}+j,\\
    \tilde{\mu}_{i-1}=\mu_i&\quad\forall i\in n\mathbb{Z}+(j-1),\\
    \tilde{\mu}_{i}=\mu_{i-1}&\quad\forall i\in n\mathbb{Z}+j,\\
    \tilde{\mu}_i=\mu_i&\quad\forall i\notin n\mathbb{Z}+(j-1)\cup n\mathbb{Z}+j.
\end{cases}\]
We write $P \overset{\mu}{\sim}_j \tilde{P}$. 

Define the group of flips to be $\Gamma_n:=\langle \phi_1,\dots,\phi_n \rangle$. Its dynamic is called the \textbf{staircase cross-ratio dynamics}.
\end{definition}

\begin{remark}
    The dynamic remains the same if we consider the discrete curvature $\lambda\mu:=(\lambda\mu_i)_{i\in\mathbb{Z}}$ for $\lambda\in\mathbb{G}_m$, since we're only interested in ratios.
\end{remark}

\begin{remark}
    This dynamical system consists of local moves, contrarily to the flat cross-ratio dynamic. Because of this, it is also well-defined if we add a monodromy $\gamma \in \mathbb{G}_m$ on the discrete curvature $\mu$:
    \[\mu_{i+n}=\gamma\mu_i, \quad \forall i\in\mathbb{Z}.\]
    Moreover, one could consider the case of \textbf{discrete curves}: these are polygons $(p_i)_{i\in\mathbb{Z}}$ endowed with discrete curvature $(\mu_i)_{i\in\mathbb{Z}}$ which are not satisfying any monodromy conditions. We won't study them here, but it might be an interesting subject. 
\end{remark}

\begin{figure}
    \centering
    \begin{tikzpicture}
        \pic[braid/width=2cm,red]{braid={a_1}};
        \pic[braid/width=2cm,red]{braid={a_1^{-1}}};
        \draw[dotted] (-4,0) -- (6,0);
        \draw[dotted] (-4,-1.5) -- (6,-1.5);
        \foreach \i in {0,-1.5}{
            \foreach \j in {-3,-1,1,3,5}{
                \filldraw[black] (\j,\i) circle(.8pt);
                }
            }
        \draw[red] (-2,0) -- (-2,-1.5);
        \draw[red] (4,0) -- (4,-1.5);

        \node at (-5,-0.75)  {$\dots$};
        \node at (7,-0.75)  {$\dots$};
        
        \node[above] at (-3,0) {$p_{j-2}$};
        \node[above] at (-1,0) {$p_{j-1}$};
        \node[above] at (1,0) {$\tilde{p}_j$};
        \node[above] at (3,0) {$p_{j+1}$};
        \node[above] at (5,0) {$p_{j+2}$};

        \node[above,red] at (-2,0) {$\mu_{j-2}$};
        \node[above,red] at (-0,0) {$\mu_{j}$};
        \node[above,red] at (2,0) {$\mu_{j-1}$};
        \node[above,red] at (4,0) {$\mu_{j+1}$};

        \node[below,red] at (-2,-1.5) {$\mu_{j-2}$};
        \node[below,red] at (-0,-1.5) {$\mu_{j-1}$};
        \node[below,red] at (2,-1.5) {$\mu_{j}$};
        \node[below,red] at (4,-1.5) {$\mu_{j+1}$};
        
        \node[below] at (-3,-1.5) {$p_{j-2}$};
        \node[below] at (-1,-1.5) {$p_{j-1}$};
        \node[below] at (1,-1.5) {$p_j$};
        \node[below] at (3,-1.5) {$p_{j+1}$};
        \node[below] at (5,-1.5) {$p_{j+2}$};
    \end{tikzpicture}
    \caption{Visualisation of the action of the flip $\phi_j$ at index $j$, read from down to up. It it similar to a braid on a cylinder, but it is important to keep in mind that it also acts on the vertices; indeed, the point $p_j$ is changed in $\tilde{p}_j=h^{\mu_j/\mu_{j-1}}_{p_{j+1},p_{j-1}}(p_{j})$.}
    \label{fig:braid}
\end{figure}

\begin{theorem}\label{thm:affine_sym_grp}
    Let $\Gamma_n=\langle \phi_1,\dots,\phi_n \rangle$ be the group generated by flips. 
    It has presentation 
    \[\langle \phi_1,\dots,\phi_n \mid \phi_i^2=1, \quad \phi_i \phi_j = \phi_j \phi_i \text{ if } |i-j|\geq2, \quad \phi_{i}\phi_{i+1}\phi_{i}=\phi_{i+1}\phi_{i}\phi_{i+1} \rangle,\]
    where the indices (and the distance) are understood modulo $n$. This is usually called the affine symmetric group (see \cite{lewisAffineSymmetricGroup2021}). 
\end{theorem}

\begin{proof}
    As it was done by Adler in \cite{adlerRecuttingsPolygons1993} for the polygon recutting, we will check each relations.
    \begin{itemize}
        \item \underline{$\phi_i^2=1$:} The action of $\phi_i$ switches $\mu_{i-1}$ with $\mu_i$, and applies $h^{\mu_i/\mu_{i-1}}_{p_{i+1},p_{i-1}}$ on $p_i$. Applying it again switches back $\mu_{i}$ with $\mu_{i-1}$, and apply $h^{\mu_{i-1}/\mu_{i}}_{p_{i+1},p_{i-1}}$ on $h^{\mu_{i}/\mu_{i-1}}_{p_{i+1},p_{i-1}}(p_i)$. So the $\mu_{i-1}, \mu_i$ are back at their place, and with the proposition \ref{prop:harm}(3) we have $h^{\mu_{i-1}/\mu_{i}}_{p_{i+1},p_{i-1}}\circ h^{\mu_{i}/\mu_{i-1}}_{p_{i+1},p_{i-1}}(p_i)=p_i$. Hence the result.
        \item \underline{$\phi_i \phi_j = \phi_j \phi_i \text{ if } |i-j|\geq2$:} This comes directly from the fact that $\phi_i$ only acts on $\mu_{i-1},\mu_i$ and $p_i$ (using $p_{i-1}, p_{i+1}$). So any change on indices $j$ with $|i-j|\geq 2$ doesn't affect anything on how $\phi_i$ will act, and we have the commutation.
        \item \underline{$\phi_{i}\phi_{i+1}\phi_{i}=\phi_{i+1}\phi_{i}\phi_{i+1}$:} This comes from the reinterpretation of a property called ``3D-consistency'' (\cite[Prop. 9]{bobenko_discrete_2008}). It states that if the vertices $(p_{i-1},p_{i},p_{i+1},p_{i+2})$ and edges $(\mu_{i-1},\mu_{i},\mu_{i+1})$ are positioned on a cube like on figure~\ref{fig:3D_consist}, one can ``flip around the cube'' by doing $(\phi_{i+1}\phi_i)^3$ and get back to the same state. So $(\phi_{i+1}\phi_i)^3=\mathrm{Id}$, and since flips are involutions we get the result.\qedhere
    \end{itemize}
\end{proof}

\begin{figure}
    \centering
    \begin{tikzpicture}
        \filldraw[red] (0,0) circle(2pt);
        \filldraw[red] (0,3) circle(2pt);
        \filldraw[red] (3,0) circle(2pt);
        \filldraw[black] (3,3) circle(2pt);
        \filldraw[black] (0+1,0+1) circle(2pt);
        \filldraw[black] (0+1,3+1) circle(2pt);
        \filldraw[red] (3+1,0+1) circle(2pt);
        \filldraw[black] (3+1,3+1) circle(2pt);

        \draw[red] (0,0) -- (0,3);
        \draw[red] (0,0) -- (3,0);
        \draw[dashed] (0,0) -- (0+1,0+1);
        \draw[dashed] (0+1,0+1) -- (0+1,3+1);
        \draw[dashed] (0+1,0+1) -- (3+1,0+1);
        \draw[red] (3,0) -- (3+1,0+1);
        \draw (3,0) -- (3,3);
        \draw (3+1,0+1) -- (3+1,3+1);
        \draw (0,3) -- (3,3);
        \draw (0+1,3+1) -- (3+1,3+1);
        \draw (0,3) -- (0+1,3+1);
        \draw (3,3) -- (3+1,3+1);

        \node[left] at (0,0) {$p_i$};
        \node[left] at (0,3) {$p_{i-1}$};
        \node[right] at (3,0) {$p_{i+1}$};
        \node[right] at (3+1,0+1) {$p_{i+2}$};
        \node[right] at (3,3) {$\tilde{p}_{i}$};
        \node[right] at (3+1,3+1) {$\tilde{p}_{i+1}$};
        \node[left] at (0+1,3+1) {$\tilde{\tilde{p}}_{i}$};
        \node[left] at (0+1,0+1) {$\tilde{\tilde{p}}_{i+1}$};

        \node[below] at (1.5,0) {$\mu_{i}$};
        \node[below] at (1.5+1,0+1) {$\mu_{i}$};
        \node[above] at (1.5,3) {$\mu_{i}$};
        \node[above] at (1.5+1,3+1) {$\mu_{i}$};

        \node[right] at (0.5,0.5) {$\mu_{i+1}$};
        \node[right] at (3.5,0.5) {$\mu_{i+1}$};
        \node[left] at (0.5,3.5) {$\mu_{i+1}$};
        \node[left] at (3.5,3.5) {$\mu_{i+1}$};

        \node[left] at (0,1.5) {$\mu_{i-1}$};
        \node[right] at (3,1.5) {$\mu_{i-1}$};
        \node[left] at (0+1,1.25+1) {$\mu_{i-1}$};
        \node[right] at (3+1,1.25+1) {$\mu_{i-1}$};
    \end{tikzpicture}
    \caption{A visualization of 3D consistency. This means that starting from the red vertices and edges, doing the sequences of flips $\phi_i\phi_{i+1}\phi_i$ or $\phi_{i+1}\phi_i\phi_{i+1}$ give the same result.}
    \label{fig:3D_consist}
\end{figure}
\begin{remark}
    The group $\Gamma_n$ doesn't always act faithfully, for instance when the field $F$ is finite (since any orbit is finite). However, when there is an element $q\in F$ of infinite multiplicative order\footnote{For instance $2\in\mathbb{Q}$.}, then each iterate of $\psi :=\phi_n\circ\dots\circ\phi_1 \in \Gamma_n$ acts differently (see subsection \ref{subsub:special}). Hence the dynamic doesn't consist only of finite orbits.
\end{remark}

Next lemma is the analogue of lemma~\ref{lemma:flat_cr}. It characterises the dynamic on the moduli space, for any flip $\phi_j$. Notice that it holds for discrete curves.

\begin{lemma}\label{lemma:stair_cr}
    Let $P=((p_1,\dots,p_n),M)$ be an $n$-gon with discrete curvature $\mu=(\mu_1,\dots,\mu_n)$. Let $j\in\mathbb{Z}$ and $\tilde{P}=\phi_j(P)$. For all $i\in \mathbb{Z}$, set
        \[c_i=[p_i, p_{i+1}, p_{i-1}, p_{i+2}],\quad \tilde{c}_i=[\tilde{p}_i, \tilde{p}_{i+1}, \tilde{p}_{i-1}, \tilde{p}_{i+2}],\]
    and set 
    \[x_j=[p_{j-2},p_{j-1},\tilde{p}_j,p_j],\quad \tilde{x}_j=[\tilde{p}_{j-2},\tilde{p}_{j-1},p_j,\tilde{p}_j],\]
    \[y_j=[p_{j+1},p_{j+2},p_j,\tilde{p}_j],\quad \tilde{y}_j=[\tilde{p}_{j+1},\tilde{p}_{j+2},\tilde{p}_j,p_j],\]
    such that 
    \[x_j \tilde{x}_j=1=y_j \tilde{y}_j.\]
    Then we have 
    \begin{align*}
        c_i&=\tilde{c}_i \qquad \quad \forall i\neq j-2,j-1,j,j+1,\tag{$\star_i$}\label{eq:stair_cr_i}\\
        c_{j-2}&=\tilde{c}_{j-2} x_j, \tag{$\star_{j-2}$}\label{eq:stair_cr_j-2}\\
        c_{j-1}(\tilde{c}_{j-1}-1)\mu_{j-1}&=\tilde{c}_{j-1}(c_{j-1}-1)\tilde{\mu}_{j-1}, \tag{$\star_{j-1}$}\label{eq:stair_cr_j-1}\\
        c_{j}(\tilde{c}_{j}-1)\mu_{j}&=\tilde{c}_{j}(c_{j}-1)\tilde{\mu}_{j}, \tag{$\star_{j}$}\label{eq:stair_cr_j}\\
        c_{j+1}&=\tilde{c}_{j+1} y_j. \tag{$\star_{j+1}$}\label{eq:stair_cr_j+1}
    \end{align*}
    Conversely, for any fixed $j\in\mathbb{Z}$, $x_j,y_j \in \mathbb{G}_m$ and $(\mu_i) \in (\mathbb{G}_m)^n$, any $(c_i)$ and $(\tilde{c}_i)$ verifying the previous equations define a pair $P \overset{\mu}{\sim}_j \tilde{P}$, unique up to a simultaneous projective transformation.
\end{lemma}

\begin{proof}
    This proof is quite similar to the one of lemma~\ref{lemma:flat_cr}. Since the flip $\phi_j$ leaves the vertices $p_i$ untouched for $i \notin n\mathbb{Z}+j$, we have that 
    \[c_i=\tilde{c}_i\qquad \forall i\neq j-2,j-1,j,j+1.\]
    Hence \ref{eq:stair_cr_i} is proved. We also obtain that $x_j=\tilde{x}_j^{-1}$ and $y_j=\tilde{y}_j^{-1}$.
    
    By the ``Chasles' relation'' of the cross ratio, we have:
    \begin{align*}
        c_{j-2}&=[p_{j-2},p_{j-1},p_{j-3},p_{j}]\\
        &=[p_{j-2},p_{j-1},p_{j-3},\tilde{p}_{j}][p_{j-2},p_{j-1},\tilde{p}_{j},p_j]\\
        &=\tilde{c}_{j-2}x_j,\\
        \frac{\mu_{j-1}}{\tilde{\mu}_{j-1}}&=[\tilde{p}_j,p_j,p_{j+1},p_{j-1}]\\
        &=[\tilde{p}_j,p_{j-2},p_{j+1},p_{j-1}][p_{j-2},p_j,p_{j+1},p_{j-1}]\\
        &=\frac{\tilde{c}_{j-1}}{\tilde{c}_{j-1}-1}\frac{c_{j-1}-1}{c_{j-1}},\\
        \frac{\mu_j}{\tilde{\mu}_{j}}&=[p_j,\tilde{p}_j,p_{j+1},p_{j-1}]\\
        &=[p_j,p_{j+2},p_{j+1},p_{j-1}][p_{j+2},\tilde{p}_j,p_{j+1},p_{j-1}]\\
        &=\frac{c_{j}-1}{c_{j}}\frac{\tilde{c}_{j}}{\tilde{c}_{j}-1},\\
        c_{j+1}&=[p_{j+1},p_{j+2},p_{j},p_{j+3}]\\
        &=[p_{j+1},p_{j+2},p_{j},\tilde{p}_{j}][p_{j+1},p_{j+2},\tilde{p}_j,p_{j+3}]\\
        &=y_j \tilde{c}_{j+1}.
    \end{align*}
    Hence we get \ref{eq:stair_cr_j-2}, \ref{eq:stair_cr_j-1}, \ref{eq:stair_cr_j}, \ref{eq:stair_cr_j+1}.\\
    
    Conversely, choose some representatives $P=(p_i)_{i\in\mathbb{Z}}$ and $\tilde{P}=(\tilde{p}_i)_{i\in\mathbb{Z}}$, equipped with discrete curvatures $\mu$ and $\tilde{\mu}$, such that $c_i=[p_i,p_{i+1},p_{i-1},p_{i+2}]$ and $\tilde{c}_i=[\tilde{p}_i,\tilde{p}_{i+1},\tilde{p}_{i-1},\tilde{p}_{i+2}]$. Because of lemma~\ref{lemma:monodromy}, we know their monodromies $M_P$, $M_{\tilde{P}}$. 
    
    The equations \ref{eq:stair_cr_i} give us equalities between many pairs of cross-ratios. Hence there exists a projective transformation $g$ such that:
    \[\tilde{p}_i=g(p_i)\quad \forall i \in\{j+1,\dots,j+n-1\}.\]
    We have:
    \begin{align*}
        c_{j+1}
        &=[p_{j+1},p_{j+2},p_{j},p_{j+3}]\\
        &=[p_{j+1},p_{j+2},p_{j},\tilde{p}_{j}][p_{j+1},p_{j+2},\tilde{p}_{j},p_{j+3}]\\
        &=y_j[g(p_{j+1}),g(p_{j+2}),g(p_{j}),g(\tilde{p}_{j+3})]\\
        &=y_j[\tilde{p}_{j+1},\tilde{p}_{j+2},\tilde{p}_{j},g(\tilde{p}_{j+3})].
    \end{align*}
    By \ref{eq:stair_cr_j+1} we obtain that 
    \[[\tilde{p}_{j+1},\tilde{p}_{j+2},\tilde{p}_{j},g(\tilde{p}_{j+3})]=\tilde{c}_{j+1}= [\tilde{p}_{j+1},\tilde{p}_{j+2},\tilde{p}_{j},\tilde{p}_{j+3}],\]
    which forces $g$ to be the identity, so
    \[\tilde{p}_i=p_i\quad \forall i \in\{j+1,\dots,j+n-1\}.\]
    By doing almost the same and using \ref{eq:stair_cr_j-2} instead, we obtain
    \[\tilde{p}_i=p_i\quad \forall i \in\{j-n+1,\dots,j-1\}.\]
    This gives us $M_P=M_{\tilde{P}}$, because they both map $(p_{j-n+1},p_{j-n+2},p_{j-n+3})$ to $(p_{j+1},p_{j+2},p_{j+3})$. Hence we have
    \[p_i=\tilde{p}_i\quad \forall i\notin n\mathbb{Z}+j.\]
    By rearranging the equation \ref{eq:stair_cr_j} in the same way as in the beginning of this proof, we see that
    \[\tilde{p}_j=h^{\mu_j/\mu_{j-1}}_{p_{j+1},p_{j-1}}(p_{j}).\]
    Hence $P \overset{\mu}{\sim}_j \tilde{P}$, which concludes the proof.
\end{proof}

Now we can use the previous lemma to get a family of scaling symmetries. They don't depend on the index $j$ where the flip is performed. Notice that again, they still hold for discrete curves. However, they are not adapted to $n$-gons with a monodromy on $\mu$.

\begin{proposition}\label{prop:stair_scaling}
    For each $\eta\in\mathbb{G}_m$, there is a scaling symmetry $\rho^{(\eta)}$ that acts in the following way:
    \[\mu_i^{(\eta)}(t)=\frac{t\mu_i}{1+(t-1)\eta\mu_i},\quad c_i^{(\eta)}(t)=c_i\frac{t \mu_i}{\mu_i^{(\eta)}(t)}.\]
    This is well defined for $t\neq1-(\eta\mu_i)^{-1}$.
\end{proposition}

\begin{remark}
    Considering the scaling symmetries for each $\eta \in \mathbb{G}_m$ is useful. Indeed, for any $\lambda\in \mathbb{G}_m$, the discrete curvatures $\mu$ and $\lambda\mu$ provide the same dynamic. They are related by this commutative diagram.
    \[\begin{tikzcd}
     (c,\mu) \arrow{r}{\rho^{(\eta)} (t)} \arrow[swap]{d}{\lambda} & \rho^{(\eta)}(t)\cdot(c,\mu) \arrow{d}{\lambda} \\
    (c,\lambda\mu) \arrow{r}{\rho^{(\eta/\lambda)} (t)} & \rho^{(\eta/\lambda)}(t)\cdot(c,\lambda\mu)
    \end{tikzcd}\]
    We will see later that the infinitesimal monodromy is independent of the choice of $\eta$.
\end{remark}

\begin{proof}
    Like for the proof of proposition \ref{prop:flat_scaling}, we have first to check that $\rho^{(\eta)}$ is indeed a group action. First, we clearly have $\mu_i^{(\eta)}(1)=\mu_i$ and $c_i^{(\eta)}(1)=c_i$. Then we have
    \begin{align*}
        \mu_i^{(\eta)}(t_1)(t_2)&=\frac{t_2\mu_i^{(\eta)}(t_1)}{1+(t_2-1)\eta\mu_i^{(\eta)}(t_1)}\\
        &=\frac{t_2t_1\mu_i}{1+(t_1-1)\eta\mu_i}\frac{1}{1+(t_2-1)\frac{t_1\mu_i}{1+(t_1-1)\eta\mu_i}}\\
        &=\frac{t_1 t_2 \mu_i}{1+(t_1-1)\eta\mu_i+t_1(t_2-1)\eta\mu_i}\\
        &=\frac{t_1t_2\mu_i}{1+(t_1t_2-1)\eta\mu_i}\\
        &=\mu_i^{(\eta)}(t_1t_2),
    \end{align*}
    and from this we get 
    \begin{align*}
        c_i^{(\eta)}(t_1)(t_2)&=c_i^{(\eta)}(t_1)\frac{t_2\mu_i^{(\eta)}(t_1)}{\mu_i^{(\eta)}(t_1)(t_2)}\\
        &=c_i\frac{t_1\mu_i}{\mu_i^{(\eta)}(t_1)}\frac{t_2\mu_i^{(\eta)}(t_1)}{\mu_i^{(\eta)}(t_1t_2)}\\
        &=c_i\frac{t_1t_2\mu_i}{\mu_i^{(\eta)}(t_1t_2)}\\
        &=c_i^{(\eta)}(t_1t_2).
    \end{align*}
    So $\rho^{(\eta)}$ acts indeed as a group. It is well defined when every $\mu_i^{(\eta)}(t)$ is in $\mathbb{G}_m$, which means that $t\neq1-(\eta\mu_i)^{-1}$.
    
    Let's show that $y_j^{(\eta)}(t)$ and $x_j^{(\eta)}(t)$ are constant. Since both cases are similar, we only perform the computation for $y_j$. To do so, take the lift\footnote{See remark \ref{remark:P_t} for details.} such that for any $t$ 
    \[(p_{j-1}^{(\eta)}(t),p_{j}^{(\eta)}(t),p_{j+1}^{(\eta)}(t))=(1,\infty,0),\]
    which forces that $p_{j+2}^{(\eta)}(t)=c_j^{(\eta)}(t)$. By proposition \ref{prop:harm}(2), we have
    \[\tilde{p}_{j}^{(\eta)}(t)=\frac{\mu_{j-1}^{(\eta)}(t)}{\mu_{j-1}^{(\eta)}(t)-\mu_{j}^{(\eta)}(t)}.\]
    From this we obtain
    \begin{align*}
        y_j^{(\eta)}(t)&=[0,p_{j+2}^{(\eta)}(t),\infty,\tilde{p}_{j}^{(\eta)}(t)]\\
        &=\frac{\tilde{p}_{j}^{(\eta)}(t)-p_{j+2}^{(\eta)}(t)}{\tilde{p}_{j}^{(\eta)}(t)}\\
        &=1-c_j\frac{t\mu_{j}}{\mu_{j}^{(\eta)}(t)}\frac{\mu_{j-1}^{(\eta)}(t)-\mu_{j}^{(\eta)}(t)}{\mu_{j-1}^{(\eta)}(t)}\\
        &=1-c_j\mu_j(\frac{1-(t-1)\eta\mu_j}{\mu_j}-\frac{1-(t-1)\eta\mu_{j-1}}{\mu_{j-1}})\\
        &=1-c_j\mu_j(\frac{1}{\mu_j}-\frac{1}{\mu_{j-1}}),
    \end{align*}
    so it is constant. Then, because of the converse part of lemma~\ref{lemma:stair_cr}, we need to check that equations \ref{eq:stair_cr_i}, \ref{eq:stair_cr_j-2}, \ref{eq:stair_cr_j-1}, \ref{eq:stair_cr_j} and \ref{eq:stair_cr_j+1} hold for the pairs $(c^{(\eta)}(t),\mu^{(\eta)}(t))$ and $(\tilde{c}^{(\eta)}(t),\tilde{\mu}^{(\eta)}(t))$.\\
    First of all, we have that $c_i^{(\eta)}(t)=\tilde{c}_i^{(\eta)}(t)$ for each $i\neq j-2,j-1,j,j+1$, so \ref{eq:stair_cr_i} is still verified.
    
    Then, we have
    \[x_i\tilde{c}_{j-2}^{(\eta)}(t)=x_i\tilde{c}_{j-2}\frac{t\mu_{j-2}}{\mu_{j-2}^{(\eta)}(t)}=c_{j-2}\frac{t\mu_{j-2}}{\mu_{j-2}^{(\eta)}(t)}=c_{j-2}^{(\eta)}(t),\]
    so \ref{eq:stair_cr_j-2} is still verified. We get \ref{eq:stair_cr_j+1} in the exact same way. Furthermore, we have
    \begin{align*}
        c_j^{(\eta)}(t)(\tilde{c}_j^{(\eta)}(t)-1)\mu_j^{(\eta)}(t)&=c_j\frac{t\mu_j}{\mu_j^{(\eta)}(t)}(\tilde{c}_j\frac{t\tilde{\mu_j}}{\tilde{\mu}_j^{(\eta)}(t)}-1)\mu_j^{(\eta)}(t)\\
        &=c_jt\mu_j(\tilde{c_j}(1+(t-1)\eta\tilde{\mu}_j)-1)\\
        &=tc_j(\tilde{c}_j-1)\mu_j + t(t-1)c_j\tilde{c}_j\eta\mu_j\tilde{\mu}_j\\
        &=t\tilde{c}_j(c_j-1)\tilde{\mu}_j + t(t-1)\tilde{c}_jc_j\eta\tilde{\mu}_j\mu_j\\
        &=\tilde{c}_jt\tilde{\mu}_j(c_j(1+(t-1)\eta\mu_j)-1)\\
        &=\tilde{c}_j^{(\eta)}(t)(c_j^{(\eta)}(t)-1)\tilde{\mu}_j^{(\eta)}(t),
    \end{align*}
    so \ref{eq:stair_cr_j} is still verified, and the same computation also gives us \ref{eq:stair_cr_j-1}.
\end{proof}

\begin{figure}
    \centering
    \includegraphics[width=0.5\linewidth]{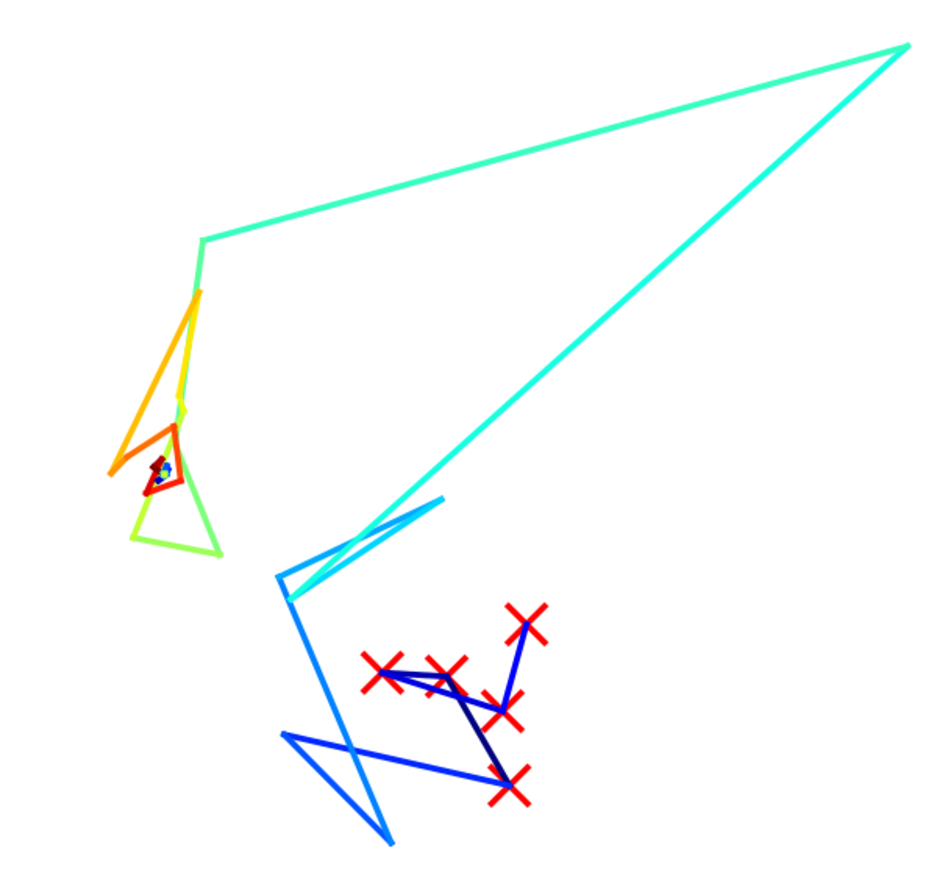}
    \caption{Python simulation of the evolution of a random closed $5$-gon in $\mathbb{P}^1(\mathbb{C})$ with random discrete curvature, under the iteration of the sequence of flips $\phi_5\circ\dots\circ\phi_1$ (applied 25 times). The initial vertices are represented with red crosses. A new edge is drawn each time a flip is done, and the colour indicates time (from blue to red). There is a rapid spiraling towards one of the collapse points determined by theorem~\ref{thm:collapse_closed}.}
    \label{fig:simulation_5-gon}
\end{figure}

Now we can compute the infinitesimal monodromy and the coordinate of the collapse points\footnote{The polynomial for the collapse points was first found experimentally by Paul Melotti when working on discrete holomorphic functions.}. Figure~\ref{fig:simulation_5-gon} depicts the convergence of a random closed $5$-gon, under the iteration of some non-trivial element of $\Gamma_5$, towards one of the collapse points, as predicted by conjecture~\ref{conj:asympt}.

\begin{proposition}
    The infinitesimal monodromy for the staircase cross-ratio dynamics, under the scaling symmetry $\rho^{(\eta)}$, is
    \[M'(1)=\sum_{i=1}^n \eta\mu_i\frac{1}{p_i-p_{i+1}}\begin{pmatrix}
        p_i & -p_ip_{i+1}\\ 
        1 & -p_{i+1}
    \end{pmatrix},\]
    and its eigenvectors $(X,Y)^T$ are roots of the polynomial 
    \[\chi_P(X,Y)=\sum_{i=1}^n \mu_i\frac{(X-p_iY)(X-p_{i+1}Y)}{p_i-p_{i+1}}.\]
\end{proposition}

\begin{remark}
    If we define $\alpha_i:=[p_i,p_{i+1},p_{i-1},\tilde{p}_i]$ to match the convention of the flat cross-ratio dynamics, we have by the permutation of the variables of the cross-ratio that
    \[\mu_i=[p_i,\tilde{p}_i,p_{i+1},p_{i-1}]=\frac{\alpha_i-1}{\alpha_i}.\]
    We get back the same formula as in proposition \ref{prop:flat_infinitesimal}. This was expectable, since they both arise from the theory of discrete holomorphic functions.
\end{remark}

\begin{proof}
    From proposition \ref{prop:stair_scaling} we have that 
    \[\frac{c_i'(1)}{c_i(1)}=\frac{c_i \eta\mu_i}{c_i}=\eta\mu_i,\]
    and we get the infinitesimal monodromy by theorem~\ref{thm:monod_closed}. We also get the polynomial $\chi_P$, which is of the desired form when we divide it by $\eta$.
\end{proof}

\subsection{Special configurations}\label{subsub:special}

\begin{figure}
    \centering
    \begin{tikzpicture}
        \draw[gray!30] (0,0) grid (2,3);
        \foreach \i in {0,1,...,2}{
            \foreach \j in {0,1,...,3}{
                %\pgfmathsetmacro\h{int(\i+\j)};
                \filldraw[black] (\i,\j) circle(.8pt);
            }
        }
        \node[right] at (1,2+.2) (0) {$q^2$};
        \node[right] at (1,1+.2) (0) {$q$};
        \node[right] at (1,0+.2) (0) {$1$};
        \node[right] at (0,3+.2) (0) {$q$};
        \node[right] at (0,2+.2) (0) {$1$};
        \node[right] at (2,0+.2) (0) {$q^2$};
    
        \node[right] at (1,3+.2) (0) {$q^3$};
        \node[right] at (2,1+.2) (0) {$q^3$};
        \node[right] at (2,2+.2) (0) {$q^4$};
        \node[right] at (2,3+.2) (0) {$q^5$};

        \node[right] at (0,1+.2) (0) {$q^{-1}$};
        \node[right] at (0,1+.2) (0) {$q^{-1}$};
        \node[right] at (0,0+.2) (0) {$q^{-2}$};
    
        \draw (0,2) -- (0,3);
        \draw (0,2) -- (1,2);
        \draw (1,2) -- (1,1);
        \draw (1,1) -- (1,0);
        \draw (1,0) -- (2,0);
        \foreach \i in {0,1,...,2}{
            \draw[dotted, red] (-0.5,\i+0.5) -- (2.5,\i+0.5);
            \node[left, red] at (-0.5,\i+0.5)  {$1$};
        }
        \foreach \i in {0,1,...,1}{
            \draw[dotted, blue] (\i+0.5,-0.5) -- (\i+0.5,3.5);
            \node[below, blue] at (\i+0.5,-0.5)  {$\lambda$};
        }

        \draw[gray!30] (0-5,0) grid (2-5,3);
        \foreach \i in {0,1,...,2}{
            \foreach \j in {0,1,...,3}{
                %\pgfmathsetmacro\h{int(\i+\j)};
                \filldraw[black] (\i-5,\j) circle(.8pt);
            }
        }
        \node[right] at (1-5,2+.2) (0) {$2$};
        \node[right] at (1-5,1+.2) (0) {$1$};
        \node[right] at (1-5,0+.2) (0) {$0$};
        \node[right] at (0-5,3+.2) (0) {$1$};
        \node[right] at (0-5,2+.2) (0) {$0$};
        \node[right] at (2-5,0+.2) (0) {$2$};
    
        \node[right] at (1-5,3+.2) (0) {$3$};
        \node[right] at (2-5,1+.2) (0) {$3$};
        \node[right] at (2-5,2+.2) (0) {$4$};
        \node[right] at (2-5,3+.2) (0) {$5$};

        \node[right] at (0-5,1+.2) (0) {$-1$};
        \node[right] at (0-5,1+.2) (0) {$-1$};
        \node[right] at (0-5,0+.2) (0) {$-2$};
    
        \draw (0-5,2) -- (0-5,3);
        \draw (0-5,2) -- (1-5,2);
        \draw (1-5,2) -- (1-5,1);
        \draw (1-5,1) -- (1-5,0);
        \draw (1-5,0) -- (2-5,0);
        \foreach \i in {0,1,...,2}{
            \draw[dotted, red] (-0.5-5,\i+0.5) -- (2.5-5,\i+0.5);
            \node[left, red] at (-0.5-5,\i+0.5)  {$1$};
        }
        \foreach \i in {0,1,...,1}{
            \draw[dotted, blue] (\i+0.5-5,-0.5) -- (\i+0.5-5,3.5);
            \node[below, blue] at (\i+0.5-5,-0.5)  {$\lambda$};
        }
        \end{tikzpicture} 
    \caption{Special ``staircase'' of the cross-ratio dynamic, for $n=3$. On the left, this makes appear a parabolic Möbius transformation, and on the right an elliptic or loxodromic (depending on $|q|$).}
    \label{fig:cliff}
\end{figure}

Here, we investigate some specific cases of staircase cross-ratio dynamics for which we can write explicitly the dynamic. Their motion on the moduli space is periodic, meaning that it amounts on the initial space to iterate a projective transformation (of any kind). They can be seen as a special case of discrete holomorphic functions, depicted on figure~\ref{fig:cliff}. The idea comes from \cite{leitenbergerIteratedHarmonicDivisions2016a}, and it holds on any algebraicly closed field of characteristic $0$. For an example of generic orbit in $\mathbb{P}^1(\mathbb{C)}$, see the simulation on figure~\ref{fig:simulation_5-gon}.\\

First, let's see that how a parabolic transformation can appear. Let's consider the closed $n$-gon 
\[P=(0,1,\dots,n-1)\]  with discrete curvature $(1,\dots,1,\lambda)$, such that after applying $\phi_n\circ\dots\circ\phi_1$, we get \[\tilde{P}=(n,n+1,\dots,2n-1).\] This happens if the cross-ratio on any square is:
\[\lambda=\frac{(k-(k+n-1)(k+n-(k+1))}{(k+n-1-(k+n))(k+1-k)}=(n-1)^2.\]
Hence, we get the Möbius transformation $z \mapsto z+n$, which is parabolic. The dynamic converges in the future and past towards the double point $\infty$.\\

Now, let's do the same for the closed $n$-gon \[P=(1,q,\dots,q^{n-1})\] with discrete curvature $(1,\dots,1,\lambda)$, such that after applying $\phi_n\circ\dots\circ\phi_1$ we get \[\tilde{P}=(q^n,q^{n+1},\dots,q^{2n-1}).\] Again, we have to check the cross-ratio condition for every square: 
\[\lambda=\frac{(q^k-q^{k+n-1})(q^{k+n}-q^{k+1})}{(q^{k+n-1}-q^{k+n})(q^{k+1}-q^k)}=\frac{(1-q^{n-1})(q^n-q)}{(q^{n-1}-q^{n})(q-1)}.\]
By remembering that $\frac{q^{n-1}-1}{q-1}=\sum_{k=0}^{n-2} q^k$ and doing some simplifications, we get that $q$ is a root of the polynomial equation
\[Q_{n}(q,\lambda)=(\sum_{k=0}^{n-2} q^k)^2-\lambda q^{n-2}=0.\]
The polynomial $Q_{n}$ is of degree $2n-4$ and palindromic with respect to $q$, which implies that $1/q$ is also a root. So we can suppose $|q|\geq1$. Note that if $q=1$, we get back $\lambda=(n-1)^2$ as in the parabolic case.\\

Let's go back to theorem \ref{thm:asympt_periodic}. In the case of periodic motion of the moduli space, any type of Möbius transformation can appear.
\begin{enumerate}
    \item If $|q|\neq1$, we get a convergence towards $\infty$ in the future and $0$ in the past (loxodromic transformation, the generic case);
    \item if $q\neq 1$ is a root of unity, then the dynamic is periodic (finite order elliptic transformation);
    \item if $|q|=1$ but is not a root of unity, then the dynamic is recurrent (infinite order elliptic transformation);
    \item otherwise, we get a convergence in the past and future towards $\infty$ (parabolic transformation).
\end{enumerate}
In all of this cases, we can deform the polygons via the scaling symmetry to get families of twisted $n$-gons whose vertices also collapse.

\appendix
\afterpage{%
    \clearpage% Flush earlier floats (otherwise order might not be correct)
    \thispagestyle{empty}% empty page style (?)
    \begin{landscape}% Landscape page
        \phantomsection
        \addcontentsline{toc}{section}{\protect\numberline{\thesection}Appendix: examples of polygonal dynamics}
        \centering % Center table
        \begin{tabular}{|p{40mm}|p{5mm}|p{15mm}|p{40mm}|p{18mm}|p{38mm}|p{60mm}|}
            \hline
            \textbf{Name of the dynamic} & $\boldsymbol X$ & $\boldsymbol G$ & \textbf{Integrability} & \textbf{Scaling symmetry} & \textbf{Group} & \textbf{Remark} \\\hline
            Leapfrog map & $\mathbb{P}^1$ & $\PGL{2}$ & Hamiltonian (\cite[§4.5]{gekhtmanIntegrableClusterDynamics2016}) & Prop.~\ref{prop:leapfrog_scaling}, adapted from \cite{gekhtmanIntegrableClusterDynamics2016} & Single transformation & Is a part of the pentagram maps family \cite[§5.2]{gekhtmanIntegrableClusterDynamics2016}\\ \hline
            Flat cross-ratio dynamic & $\mathbb{P}^1$ & $\PGL{2}$ & Hamiltonian (\cite{arnoldCrossratioDynamicsIdeal2022},\cite{affolterIntegrableDynamicsProjective2023}) and algebraic (\cite{hetrich-jerominPeriodicDiscreteConformal2001})& Prop.~\ref{prop:flat_scaling}& Single transformation & Has a discrete curvature\\ \hline
            Staircase cross-ratio dynamic & $\mathbb{P}^1$ & $\PGL{2}$ & Probably algebraic (for a future paper)& Prop.~\ref{prop:stair_scaling} & Affine symmetric group (Thm.~\ref{thm:affine_sym_grp}) & Has a discrete curvature, first defined here\\ \hline
            Pentagram map & $\mathbb{P}^2$ & $\PGL{3}$ & Hamiltonian (\cite{ovsienkoPentagramMapDiscrete2010},\cite{ovsienkoLiouvilleArnoldIntegrabilityPentagram2011}) and algebraic (\cite{solovievIntegrabilityPentagramMap2013},\cite{weinreichAlgebraicDynamicsPentagram2023}) & \cite[Cor. 2.5]{ovsienkoPentagramMapDiscrete2010} & Single transformation & Introduced in \cite{schwartzPentagramMap1992}, survey in \cite{stieglerwikijournal2025}\\ \hline
            Generalized pentagram maps & $\mathbb{P}^d$ & $\PGL{d+1}$ & Depends (\cite{khesin_non-integrability_2015}) & Sometimes (see \cite{beffaIntegrableGeneralizationsPentagram2015}, \cite{gekhtmanIntegrableClusterDynamics2016}) & Single transformation & See \cite{beffaIntegrableGeneralizationsPentagram2015},\cite{izosimovLongdiagonalPentagramMaps2023}, \cite{handPentagramMapsRings2024}\\ \hline
            Projective heat map & $\mathbb{P}^2$ & $\PGL{3}$ & No \cite[§15.7]{schwartzProjectiveHeatMap2017} & No & Single transformation, not birational & Introduced in \cite{schwartzProjectiveHeatMap2017}, studied in \cite{kaschnerComplexPerspectiveProjective2017}, \cite{leFamilyProjectivelyNatural2018} \\ \hline
            Polygon recutting & $\mathbb{R}^2$ & $\mathrm{Isom}^+(\mathbb{R}^2)$ & Hamiltonian (\cite{izosimovPolygonRecuttingCluster2023}) & Unknown & Affine symmetric group (see \cite{adlerRecuttingsPolygons1993})& Introduced in \cite{adlerRecuttingsPolygons1993}\\ \hline
            Polygon folding & $\mathbb{R}^2$ & $\mathrm{Isom}^+(\mathbb{R}^2)$ & No & Unknown & Folding group (see \cite[§5.4]{cantatDynamicsAutomorphismGroups2023a}) & 4-gons are treated in \cite{benoistIterationPliagesQuadrilateres2004} and there are probabilistic result for 5-gons in \cite{cantatDynamicsAutomorphismGroups2023a}\\ \hline
        \end{tabular}
        \captionof{table}{Examples of polygonal dynamics.}    \label{tab:examples_polygonal_dynamics}
    \end{landscape}
    \clearpage
}
\newpage

\subsection*{Funding}
This work was partially funded by the French National Research Agency under the project DynAtrois (ANR-24-CE40-1163).

\subsection*{Acknowledgements}
I warmly thank my PhD advisors Thomas Gauthier and Ramanujan Santharoubane for their constant support, and also Paul Melotti with whom I started working on this question during my master's. He was the first to have the intuition of the collapse formula for the system that then became the staircase cross-ratio dynamics. I'd also like to thank Niklas Affolter, Anton Izosimov, Sanjay Ramassamy, Richard Schwartz, Valentin Ovsienko and Max Weinreich for the helpful and sympathetic discussions.

\printbibliography
\end{document}